\documentclass[11pt,a4paper
]{article}

\usepackage{amsmath,amsfonts,amssymb,latexsym}
\usepackage{mathrsfs}
\usepackage[dvips]{graphicx}
\usepackage{enumerate}

\newtheorem{theorem}{Theorem}[section]
\newtheorem{lemma}[theorem]{Lemma}

\newtheorem{example}[theorem]{Example}
\newtheorem{definition}[theorem]{Definition}

\newtheorem{remark}[theorem]{Remark}

\def\R{\mathbb{R}}

\def\N{\mathbb{N}}
\def\ind{\mathbb{I}}
\def\Argmin{\rm{Argmin\; }}

\def\endproof{\hfill\Box}

\title{\small\bf
CONTINUITY AND MONOTONICITY OF %GREEDY MAXIMIZERS
SOLUTIONS TO A GREEDY MAXIMIZATION PROBLEM}

\author{ {\small\sc
\L ukasz Kruk} \\
{\em\footnotesize Maria Curie-Sk\l odowska
University, Lublin, Poland}
}

\mathchardef\emptyset="001F

\begin{document}

\maketitle

\begin{abstract}
Motivated by an application to resource sharing network modelling, we consider a problem of greedy maximization (i.e., maximization of the consecutive minima) of a vector in $\R^n$, with the admissible set indexed by the time parameter. The structure of the constraints depends on the underlying network topology. We investigate continuity and monotonicity of the resulting maximizers with respect to time. Our results have important consequences for fluid models of the corresponding networks which are optimal, in the appropriate sense, with respect to handling real-time transmission requests.
\end{abstract}



\section{Introduction}\label{s.intro}

Starting from seminal papers of Rybko and Stolyar \cite{rysto}, Dai \cite{dai}, Bramson \cite{bram, bram2, bram3} and other authors, fluid models have become a standard tool in investigating long-time behaviour of complicated
queueing %networks.
systems.
Such models are useful for establishing stability and %characterizing
obtaining hydrodynamic or %even
diffusion limits for multiclass queueing networks and resource sharing networks with various service protocols. Using a similar
methodology in the %important
case of real-time Earliest Deadline First (EDF) networks with resource sharing is hindered by the
%fact that it is not clear how to find
lack of suitable fluid model equations. To overcome this difficulty, Kruk \cite{lk4} suggested the %characterization
definition of fluid models for these systems by means of an optimality property, called local edge minimality, which is known to
characterize the EDF discipline in %actual
stochastic
%(pre-limit)
resource sharing networks. It turns out that the success of this approach depends %in a crucial way
on establishing suitable local properties of a vector-valued mapping $F:[0,\infty)\rightarrow \R^n$,
resulting from a greedy maximization (i.e., maximization of the consecutive minima) of a vector in $\R^n$ over the admissible set $A_t$, depending on the underlying network topology and indexed by the time parameter.
%moze cos o posetach tutaj?
A convenient way to describe the value $F(t)$ for a given $t\geq 0$ is to define it as the maximal element of %the admissible set
$A_t$ with respect to a suitable ``min-sensitive'' partial ordering.
Roughtly speaking, the function $F$, when well behaved, determines the so-called frontiers (i.e., the left endpoints of the supports) of the states in the corresponding locally edge minimal fluid model. The idea of using frontiers for the asymptotics of EDF %networks
systems dates back to the paper of
Doytchinov et al. \cite{DLS2001} on a G/G/1 queue, and it has been used several times since then. However, both %its
the application of this idea to resource sharing networks and our %algorithm for
approach to determining
%finding
the frontiers %vectors
%via
%linking it to
by finding
maximal elements of partially ordered sets
appear to be new.

In this paper, for any $t\geq 0$, we
%is devoted to the
construct the value of
%ion of the maximizer
$F(t)$ as a solution of a nested sequence of max-min problems in a $t$-dependent admissible set. While none of these max-min problems is %particularly
hard to solve, %in general,
their number and %exact
forms vary with $t$ in a complex, discontinuous %and rather
way, making the analysis of the resulting mapping $F$ on $[0,\infty)$ rather involved.
We %also
investigate key properties of $F$,
%which are important from the point of view of our queueing application,
namely, its continuity and monotonicity.
Our main results are described in more detail in Section \ref{ss.mainres}, to follow, after the introduction of %necessary
indispensable notation.
%of the resulting maximizers with respect to time.
These results are used in
%a subsequent
our forthcoming paper %\cite{lk5}
to establish fundamental properties of locally edge minimal fluid models, like their existence, uniqueness and stability.

We hope that the
%mapping $F$ investigated
theory developed in this paper
%and its modifications
will be useful not only in the asymptotic analysis of EDF-like
disciplines, but also in the case of other ``greedy``
scheduling policies for resource sharing networks, for example
Longest Queue First %(LQF)
%(see, e.g., %Dimakis and Walrand
\cite{DW} %)
or Shortest Remaining
Processing Time %(SRPT) %(see Verloop et al.
\cite{vbn}.

\subsection{Notation}\label{sec:not}

%The following notation will be used throughout.
For sets $A$, $B$, we write $A \subsetneq B$ if $A$ is a proper subset of $B$.
For a finite set $A$, let $|A|$ denote the cardinality of $A$.
% and let $2^A$ denote the family of all the subsets of $A$.
Let $\N$ denote the
set of positive integers and let $\R$ denote the set of real
numbers. For $a,b\in\R$, we write $a\vee b$ ($a\wedge b$) for the
maximum (minimum) of $a$
and $b$, $a^+$ for $a \vee 0$. %, $a^-$ for $(-a) \vee 0$
Vector inequalities are to be interpreted componentwise, i.e., for  $a,b\in \R^n$, $a=(a_1,..,a_n)$, $b=(b_1,...,b_n)$, $a\leq b$ if and only if $a_i\leq b_i$ for all $i=1,...,n$.
For $a=(a_1,..,a_n)\in \R^n$, we write $\min a$ for $\min_{i=1,...,n} a_i$ and
$\Argmin a$ for $\{i\in \{1,...,n\}: a_i = \min a\}$.
For $a=(a_1,..,a_n)\in \R^n$ and a set $A\subsetneq \{1,...,n\}$ with $ \{1,...,n\} \setminus A= \{i_1,...,i_k\}$, where $k=n-|A|$ and $i_1<i_2<...<i_k$,  we identify $(a_i)_{i\notin A}$ with $(a_{i_1},...a_{i_k})\in \R^k$.
By
convention,
a sum of the form $\sum_{i=n}^m$ ($\bigcup_{i=n}^m$) with $n>m$ or, more generally,
a sum of numbers (resp., sets)
over the empty set of indices equals zero (resp., $\emptyset$).
For a set $A \subseteq \R$, let $\overline{A}$ denote the closure of $A$.

\section{The %problem statement
mapping $F$}\label{s.ourmap}

\subsection{%``min-like''
``Min-sensitive'' partial ordering on $\R^n$}

We define the relation ``$\eqslantless$'' inductively on $\R^n$ as follows.
\begin{definition}\label{d.porder}
For $a,b\in \R$, we write $a \eqslantless b$ iff $a\leq b$.
If $n\geq 2$, for $a,b\in \R^n$, $a=(a_1,..a_n)$, $b=(b_1,...,b_n)$, we write $a \eqslantless b$ if one
%in each
of the following four cases holds.

(i) $a=b$,

(ii) $\min a < \min b$,

(iii) $\min a = \min b$ and $\Argmin b \subsetneq \Argmin a$,

(iv) $\min a = \min b$, $\Argmin a = \Argmin b \subsetneq \{1,...,n\}$ and $(a_i)_{i \notin \Argmin a} \eqslantless (b_i)_{i \notin \Argmin b}$ in $\R^{|\Argmin a|}$.
\end{definition}
%In the point (iv) of Definition \ref{d.porder}, we have implicitly identified
\begin{remark}\label{r.obv1}
{\rm
%The necessary (but sufficient only in the case of $n=1$) condition for $a \eqslantless b$ is $\min a \leq \min b$.
Clearly, if $a \eqslantless b$, then $\min a \leq \min b$. The converse is, in general, false, unless $n=1$. For example, the relation $(0,1) \eqslantless (1,0)$ does not hold.
}
\end{remark}
The proof of the following lemmma is elementary and it is left to the reader.
\begin{lemma}\label{l.porder}
The relation ``$\eqslantless$'' is a partial ordering on $\R^n$.
\end{lemma}
\begin{remark}\label{r.obv2}
{\rm
For $a,b\in \R^n$, the inequality $a\leq b$ implies $a \eqslantless b$. %For $n\geq 2$,
In dimensions greater than one the converse is in general, false, for example, $(0,2) \eqslantless (1,1)$.
}
\end{remark}

\subsection{The mapping definition}\label{ss.Fpose}

Let $I,J\in \N$ and let ${\bf I}=\{1,...,I\}$, ${\bf
J}=\{1,...,J\}$. For $i\in \bf{I}$, let $h_i:\R\rightarrow \R$ be a
continuous, nonnegative, nondecreasing function with
$\lim_{x\rightarrow \infty} h_i(x)=\infty$ and
\begin{equation}\label{e.defxi}
x^*_i:=\sup\{x\in \R:h_i(x)=0\} \in (-\infty, 0].
\end{equation}
In other words, each $h_i$ is the cumulative distribution function
of an atomless, $\sigma$-finite measure in $\R$, with finite,
nonpositive infimum of its support. Let $G_j$, $j\in {\bf J}$, be a
family of distinct, nonempty subsets
of ${\bf I}$ ({\em not necessarily pairwise disjoint}) %and
such that
${\bf I}= \bigcup_{j\in {\bf J}} G_j$.
\begin{definition}\label{d.F}
For %a given
$t\geq 0$, we denote by $F(t)=(F_i(t))_{i\in {\bf I}}$ the maximal element of the set
%of vectors $a=(a_1,...,a_I)\in \R^I$ satisfying the constraints
$$
A_t = \Big\{a=(a_1,...,a_I)\in \R^I:\;  a_i \leq t, \; i\in {\bf I}, \;  \sum_{i\in G_j} h_i(a_i) \leq t, \; j \in {\bf J} \Big\},
$$
with respect to the relation ``$\eqslantless$''.
\end{definition}
Somewhat informally, %the vector
$F(t)$ may be thought of as the result of ``greedy'' maximization of a vector $a\in \R^I$, subject to the constraints defining the set $A_t$, in the following sense. We first maximize $\min a$ over $a\in A_t$, then we maximize the ``next minimum'' $\min \{ a_i, i\notin \Argmin a\}$ over the set of maximizers of the previous problem, and we continue in this way until all the entries of the maximizer $a^*=F(t)$ are determined.
In Section \ref{ss.constr} %, to follow,
we formalize and describe in detail this nested %maximization
max-min procedure
%provide the construction of this maximal element
which implies,
%implying,
in particular, the existence and uniqueness of the maximizer $F(t)$.
\begin{remark}\label{r.shift}
{\rm
A seemingly more general version of Definition \ref{d.F}, in which
for a %given,
fixed (possibly positive) $t_0\in \R$ and %any
$t\geq t_0$,
the set $A_t$ is replaced by
$$
A_t^{t_0} = \Big\{a=(a_1,...,a_I)\in \R^I:\;  a_i \leq t, \; i\in {\bf I}, \;  \sum_{i\in G_j} h_i(a_i) \leq t-t_0, \; j \in {\bf J} \Big\},
$$
and the %numbers $x^*_i$ defined by
half-line $(-\infty,0]$ in \eqref{e.defxi} %are required to belong to
is replaced by $(-\infty, t_0]$,
may be easily reduced to the case %of $t_0=0$
considered in Definition \ref{d.F} %($t_0=0$)
by the %shift
change of variables $y=x-t_0$, $s=t-t_0$ and by using $\tilde{h}_i(y)=h_i(x)=h_i(y+t_0)$ instead of $h_i$.
}

\end{remark}

\subsection{Motivation: fluid models for resource sharing networks}\label{ss.moti}

The need for investigating the properties of the mapping $F$ defined %in Section \ref{ss.Fpose}
above arises %naturally
in %modeling
the theory of fluid models for real-time networks with resource sharing.
% in which priority is given to tasks with the earliest due dates.
Below, we briefly (and somewhat informally) describe this
connection. The %interested
reader may consult %Kruk
\cite{lk4} %, \cite{lk5}
for more details and references.

Consider a network with a finite number of resources (nodes), labelled by $j=1,...,J$, and a finite set of routes, labelled by $i=1,...,I$.
%Each route
%may be identified with a nonempty subset of ${\bf J}=\{1,...,J\}$, interpreted as the set of resources used by this route.
%Let $A=[a_{ji}]$ be the $J\times I$ incidence matrix in which $a_{ji}=1$ if resource $j$ is used by route $i$ and $a_{ji}=0$ otherwise.
Let ${\bf I}=\{1,...,I\}$, ${\bf J}=\{1,...,J\}$. For $j\in {\bf
J}$, let $G_j \subseteq {\bf I}$ be the set of routes using the
resource $j$. For convenience, we assume that all the resources have
a unit service rate. By a flow on route $i$ we mean a continuous
transmission of a file through the resources used by this route. We
assume that a flow takes simultaneous possession of all the
resources on its route during the transmission. Each flow in the
network has a deadline for transmission completion. Networks of this
type may be used to model, e.g., voice and video transmission,
manufacturing systems with order due dates or emergency health care
services. In what follows, by the lead time of a flow we mean the
difference between its deadline and the current time.

%By the lead time of a flow we mean the difference between its deadline and the current time.
As in %Kruk
\cite{lk3, lk4}, the time evolution of such a system may be described by the process
$ %\begin{equation}\label{e.prelimitX}
\mathfrak{X}(t,s)=(Z(t,s),D(t,s),T(t,s),Y(t,s))$, %\qquad
$t\geq 0,  s\in \R,
$ %\end{equation}
where the component processes $Z$, $D$, $T$, $Y$ are defined as follows.
For %every
$t\geq 0$ and $s\in \R$, $Z(t,s)=(Z_i(t,s))_{i\in {\bf I}}$, where
$Z_i(t,s)$ is the number of flows on route $i$ with lead times at time $t$ less than or equal to $s$ which are still present in the system at that time.
Similarly, the vectors
$D(t,s)=(D_i(t,s))_{i\in {\bf I}}$,
$T(t,s)=(T_i(t,s))_{i\in {\bf I}}$ denote the number of departures (i.e., transmission completions) and the cumulative %service
transmission time
by time $t$ corresponding to each route $i$ of flows with lead times at time $t$
less than or equal to $s$. Let $Y_i(t,s)=t-T_i(t,s)$, $i\in {\bf I}$, denote the cumulative idleness by time $t$ with regard to
transmission of flows on route $i$ with lead times at time $t$ less than or equal to $s$ and let $Y(t,s)=(Y_i(t,s))_{i\in {\bf I}}$.
The process $\mathfrak{X}$ satisfies the following network equations, valid for $\tilde{t}\geq t\geq 0$, $s\in \R$:
\begin{eqnarray}
Z(t,s)= Z(0,t+s)+E(t,s)-D(t,s), \qquad\;\; \qquad %\qquad\qquad\qquad \;\;\;\;
\label{1}\\
D_i(t,s)=S_i(T_i(t,s),t,s), \qquad %i\in {\bf I}, \;\;\qquad\;\; \qquad %\qquad \qquad
%\label{2}\\
T_{i}(t,s) +Y_i(t,s)=t, \qquad i\in {\bf I}, %\qquad  \qquad \; \qquad %\qquad\quad \;\;\;
\label{3}\\
\sum_{i\in G_j } \big( T_i(\tilde{t}, s-\tilde{t})- T_i(t,s-t) \big) \leq \tilde{t}-t, \qquad  j\in {\bf J}, \qquad \label{3.5}
\iffalse \\
\int _0^t Z_i(u,s-u) \ind_{[\sum_{k\in \mathcal{K}(i)} Z_k(u,s-u)=0]} \; d_u Y_i(u,s-u)=0, \qquad i\in {\bf I},
 \label{4}
\fi
\end{eqnarray}
where $E(t,s)=(E_i(t,s))_{i\in {\bf I}}$ is the corresponding external arrival process and
%$i\in {\bf I}$, $t,t'\geq 0$ and $s\in \R$, let
$S_i(t',t,s)$ denotes the number of
transmission completions of flows on route $i\in {\bf I}$ having lead times at time $t$
less than or equal to $s$, by the time the system has spent $t'$ units of time transmitting these flows.

Fluid models are deterministic, continuous analogs of resource sharing networks, in which individual flows are replaced by a divisible commodity (fluid), moving along $I$ routes with $J$ resources (nodes).
They usually arise from formal functional law of large numbers approximations of the corresponding stochastic flow level models.
The %deterministic
analogs of the network equations \eqref{1}-\eqref{3.5} are the following fluid model equations, valid for $\tilde{t}\geq t\geq 0$, $s\in \R$:
\begin{eqnarray}
\overline{Z}(t,s)= \overline{Z}(0,t+s)+ \alpha (t+ (s\wedge 0))^+
% ogolniej: bylo \alpha \circ \int^t_0 G(s+\eta)d\eta
-\overline{D}(t,s), \qquad\;\;\; %\qquad\qquad\qquad\qquad \;\;\;\;
\label{f1}\\
% zrobilismy  \alpha (t+ (s\wedge 0))^+
\overline{D}_i(t,s)=\overline{T}_i(t,s)/m_i, \qquad %i\in {\bf I}, \;\;\qquad\; \quad \qquad %\qquad
%\label{f2}\\
\overline{T}_{i}(t,s) +\overline{Y}_i(t,s)=t, \qquad i\in {\bf I}, %\qquad  \quad %\;\;\; \qquad %\qquad \quad \;\;\;\;
\label{f3}\\
\sum_{i\in G_j} \big( \overline{T}_i(\tilde{t}, s-\tilde{t})- \overline{T}_i(t,s-t) \big) \leq \tilde{t}-t, \qquad  j\in {\bf J}, \qquad \;\;\label{f3.5}
\iffalse
\\
\int _0^t \overline{Z}_i(u,s-u) \ind_{[\sum_{k\in \mathcal{K}(i)} \overline{Z}_k(u,s-u)=0]} \; d_u \overline{Y}_i(u,s-u)=0, \qquad i\in {\bf I},
 \label{f4}
\fi
\end{eqnarray}
where $m_i$ is the mean transmission time of a flow on route $i$ and $\alpha=(\alpha_i)_{i\in {\bf I}}$ is the vector of flow arrival rates.
%(more general versions of \eqref{f1}, corresponding to time-varying arrival rates may also by considered).
A system
\begin{equation}\label{X}
\overline{\mathfrak{X}}(t,s)=(\overline{Z}(t,s),\overline{D}(t,s),\overline{T}(t,s),\overline{Y}(t,s)), \qquad t\geq 0, \; s\in \R,
\end{equation}
with continuous, nonnegative components,
satisfying the equations \eqref{f1}-\eqref{f3.5}, together with some natural monotonicity assumptions, is called a {\em fluid model for the %EDF
resource sharing network} under consideration.
% (or just a {\em %EDF resource sharing fluid model}).

To proceed further, we will introduce a class of %such
fluid models which is, in some sense, optimal with respect to handling real-time transmission requests.
To this end,
we define a partial ordering ``$\ll$'' on the space of real functions on $\R$, which is extremely sensitive to the behaviour of the functions under comparison for small arguments.
\begin{definition}[\cite{lk4}, Definition 5]\label{d.edge<=}
Let $f, g:\R\rightarrow \R$ be such that for some $a\in \R$ we have
$f\equiv g$ on $(-\infty,a]$ and let
$
c= %c_{f,g} =
\sup \{ a\in \R: f(x)=g(x) \; \forall x %\in (-\infty,a]
\leq a\}.
$
We write $f \ll g$ if either $c=\infty$ (i.e., $f\equiv g$ on $\R$), or
$c<\infty$ and there exists $\epsilon %=\epsilon_{f,g}
>0$ such that $f\leq g$ on
$[c,c+\epsilon]$.
\end{definition}

\begin{definition}[see \cite{lk4}, Definition 11]\label{d.locminedgefluid}
A fluid model $\overline{\mathfrak{X}}$ of the form \eqref{X} for a %n EDF
resource sharing network %satisfying
with $\sum_{i\in {\bf I}} \overline{Z}_i(0,\cdot)\equiv 0$ on $(-\infty,c]$ for some $c\in \R$  %and all $i\in {\bf I}$
is called locally edge minimal at a time
%point %if for any
$t_0\geq 0$ if %, $\omega \in \Omega$
there exists $h>0$ such that for any fluid model
$\overline{\mathfrak{X}}'$ with the same
%fluid model data
$\alpha$, $m_i$, $G_i$ %as $\overline{\mathfrak{X}}$,
and the same initial state
$\overline{Z}'(0,\cdot)=\overline{Z}(0,\cdot)$,
satisfying
$ %\begin{equation}\label{e.samestartfluid}
\overline{\mathfrak{X}}'(t_0,\cdot)=\overline{\mathfrak{X}}(t_0,\cdot),
$ %\end{equation}
we have $\sum_{i\in {\bf I}} \overline{Y}_i(t,\cdot)\ll \sum_{i\in
{\bf I}} \overline{Y}_i'(t,\cdot)$
(equivalently,
%$\sum_{i\in {\bf I}} \overline{T}_i'(t,\cdot)\ll \sum_{i\in {\bf I}}
%\overline{T}_i(t,\cdot)$, %$\sum_{i\in {\bf I}}
%\overline{D}_i'(t,\cdot)\ll \sum_{i\in {\bf I}}
%\overline{D}_i(t,\cdot)$ or
$\sum_{i\in {\bf I}}
\overline{Z}_i(t,\cdot)\ll \sum_{i\in {\bf I}}
\overline{Z}_i'(t,\cdot)$)
for every $t\in (t_0,t_0+h)$.
The fluid model $\overline{\mathfrak{X}}$ is called locally edge
minimal, if it is locally edge minimal at every $t_0\geq 0$.
\end{definition}
The intuition behind these notions is that a locally edge minimal fluid model tries to transmit as much ``customer mass'' corresponding to the earliest deadlines as possible, and hence its idleness with respect to such mass is as small as possible. Accordingly, such a model may be thought of as a ``macroscopic'' counterpart of a resource sharing network working under the Earliest-Deadline-First (EDF) protocol. Indeed, %as it was shown in %Section 5 of \cite{lk4}, under mild assumptions,
the EDF service discipline in such a network may be characterized by
an analogous notion of %its
local edge minimality, see %Kruk
\cite{lk4}, Definition 8 and Theorems 5-7.

%It is shown in Kruk \cite{lk5}
In a forthcoming paper, we show that existence and local uniqueness of a  fluid model for given data $\alpha$, $m_i$, $G_i$, %as $\overline{\mathfrak{X}}$,
and %the same
initial state
$\overline{Z}(0,\cdot)$,
which is locally edge minimal at a time $t_0\geq 0$,
is closely related to local monotonicity of the mapping $F$ introduced in Section \ref{ss.Fpose}, with
\begin{equation}\label{e.ourhi}
h_i(x)=m_i [\overline{Z}_i(0,x) - \overline{Z}_i(0,x^*_i)]^+ + \alpha_i m_i [x^+-(x^*_i)^+]^+, \qquad i\in {\bf I},
\end{equation}
and suitably defined, not necessarily nonpositive, $x^*_i$, %which
depending on %the point
$t_0$ (see Remark \ref{r.shift}).
%at which local minimality of the fluid model is considered.
%Moreover, i
It turns out that for large $t_0$, the points $x^*_i$ are also large, so under a natural assumption that the supports of
$\overline{Z}_i(0,\cdot)$ are bounded above, the functions $h_i$ in \eqref{e.ourhi} are linear on $[x^*_i, \infty)$. Consequently, the
linear case, investigated in Section \ref{ss.plin} of this paper, is of %particular
considerable importance, because it determines the long-time
behaviour of the corresponding locally edge minimal fluid model, for
example, its stability or the form of its invariant manifold. In
particular, one of remarkable implications of the formulae %from
developed in Section \ref{ss.plin} is stability of locally edge
minimal fluid models, regardless of the underlying resource sharing
network topology. Results along this line may be found in
%Kruk \cite{lk5}.
our forthcoming paper

\subsection{Overwiev of the main results}\label{ss.mainres}

In this paper, we investigate key properties of the mapping $F: [0,\infty)\rightarrow \R^I$, %such as its continuity and monotonicity,
in particular those
which are %highly
%the most
relevant to the theory of locally edge minimal fluid models.
%for resource sharing networks.
%First, i
In Section \ref{s.constra}, we present a detailed construction of $F(t)$ and we provide some illustrating examples.
The main result of Section \ref{s.cont} is Theorem \ref{t.ciag}, stating that if each function $h_i$ is strictly increasing in $[x^*_i,\infty)$, then the corresponding map $F$ is continuous. This basic regularity result, together with the method of partions, introduced in Section \ref{ss.constJ}, is useful in proving various refinements, like Lipschitz continuity of $F$ for Lipschitz $h_i$ (Section \ref{ss.Lip}) or an upgrade of the local monotonicity result from %the case of
linear to $C^1$ functions $h_i$ (Section \ref{ss.C1mon}). The main contribution of Section \ref{s.mono} is the explicit evaluation of the mapping $F$ %in the neighbourhood of
near zero in the linear case, implying, in particular, its local monotonicity in a neighbourhood of zero for piecewise linear $h_i$.
As we have already mentioned, the latter fact is then generalized to (piecewise) %locally)
$C^1$ functions.
%We conclude by providing an example of
Finally, in Section \ref{ss.notglobal}, we show that, somewhat surprisingly, the mapping $F$ may %not
fail to be globally monotone on $[0,\infty)$, %the whole nonnegative half-line,
even if the corresponding functions $h_i$ are linear in $[x^*_i,\infty)$.

%The following theorems are our main results.

\section{The mapping construction algorithm}\label{s.constra}
\setcounter{equation}{0}

\subsection{Construction}\label{ss.constr}

Fix $t\geq
0$. We define the vector $F(t)=(F_i(t))_{i\in {\bf I}}$ as follows.

Let $f^{(1)}=f^{(1)}(t)$ be the supremum of $x\leq t$ satisfying the
constraints
\begin{equation}\label{e.max1}
\sum_{i\in G_j} h_i(x) \leq t, \qquad \quad j \in {\bf J}.
\end{equation}
If $f^{(1)}=t$, we take $F_i(t)=t$ for each $i\in {\bf I}$, ${\bf
I}^{(1)}={\bf I}$, ${\bf J}^{(1)}={\bf J}$, ${\bf
N}^{(1)}=\emptyset$ and $k_{max}=k_{max}(t)=1$. In what follows, we
assume that $f^{(1)}<t$. By continuity of $h_i$, $x=f^{(1)}$
satisfies \eqref{e.max1} and the set
$
{\bf J}^{(1)} = \big\{ j\in {\bf J}: \sum_{i\in G_j} h_i(f^{(1)}) = t \big\}
$
of active constraints %in \eqref{e.max1}
is nonempty. (Indeed, if ${\bf J}^{(1)} =\emptyset$, then
$f=f^{(1)}+\epsilon$ also satisfies \eqref{e.max1} for $\epsilon>0$
small enough, which contradicts the definition of $f^{(1)}$.) Let $
{\bf I}^{(1)} = \bigcup_{j\in {\bf J}^{(1)}} G_j$ and
\begin{equation}\label{e.defn1}
{\bf
N}^{(1)}= \{ j\in {\bf J}: G_j \subseteq {\bf I}^{(1)}\} \setminus
{\bf J}^{(1)}.
\end{equation}
For $i\in {\bf I}^{(1)}$, put
$ %\begin{equation}\label{e.firstF}
F_i(t)=\sup\{ x\leq t:h_i(x)=h_i(f^{(1)})\}.
$ %\end{equation}
% Z_i(t,s)= arrival measure cut from below at f^{(1)}
If ${\bf I}^{(1)}={\bf I}$, %we have defined
this completely determines the vector $F(t)$. In this case, let
$k_{max}=1$. Otherwise, let ${\bf K}^{(1)}={\bf J}\setminus ({\bf
J}^{(1)}\cup {\bf N}^{(1)})$ and let $f^{(2)}=f^{(2)}(t)$ be the
supremum of $x\leq t$ satisfying the constraints
\begin{eqnarray}
\sum_{i\in G_j\cap {\bf I}^{(1)}}
% h_i(f^{(1)})
h_i(F_i(t))
+ \sum_{i\in
G_j\setminus {\bf I}^{(1)}} h_i(x) \leq t,
\qquad \qquad j \in {\bf K}^{(1)}.
\label{e.max2}
\end{eqnarray}
%for every $j \in {\bf K}^{(1)}$.
Note that $G_j\setminus {\bf I}^{(1)}\neq \emptyset$ for $j \in {\bf
K}^{(1)}$, so the second sum in \eqref{e.max2} is taken over a
nonempty set of indices. We also have $f^{(1)} < f^{(2)}$ by
definition. If $f^{(2)}=t$, we take ${\bf I}^{(2)} = {\bf
I}\setminus {\bf I}^{(1)}$, ${\bf J}^{(2)}={\bf K}^{(1)}$, ${\bf
N}^{(2)}=\emptyset$, $F_i(t)=t$ for each $i\in {\bf I}^{(2)}$ and
$k_{max}=2$. If $f^{(2)}<t$,
%by continuity of $g^i_t(\cdot)$,
then $x=f^{(2)}$ satisfies \eqref{e.max2} and the set ${\bf
J}^{(2)}$ of active constraints in \eqref{e.max2} (i.e., those $j
\in {\bf K}^{(1)}$, for which equality holds in \eqref{e.max2} with
$x=f^{(2)}$) is nonempty. In this case, let ${\bf I}^{(2)} =
\bigcup_{j\in {\bf J}^{(2)}} G_j \setminus {\bf I}^{(1)}$, ${\bf
N}^{(2)}= \{ j\in {\bf J}: G_j \subseteq {\bf I}^{(1)}\cup {\bf
I}^{(2)}\} \setminus ({\bf J}^{(1)}\cup {\bf J}^{(2)} \cup{\bf
N}^{(1)})$ and put
$
F_i(t)=\sup\{ x\leq t:h_i(x)=h_i(f^{(2)})\}
$
for $i \in {\bf I}^{(2)}$. If ${\bf I}^{(1)}\cup {\bf I}^{(2)}={\bf
I}$, the definition of the vector $F(t)$ is complete and we take
$k_{max}=2$, otherwise we let ${\bf K}^{(2)}={\bf J}\setminus ( {\bf
J}^{(1)} \cup {\bf J}^{(2)} \cup {\bf N}^{(1)} \cup {\bf N}^{(2)})$,
and we continue our construction as follows.

Suppose that for some $k\geq 2$, we have %already
defined numbers $f^{(1)}<f^{(2)}<...<f^{(k)}<t$, 
% poprawione, bylo 0 w submitted version !!!
nonempty, disjoint subsets ${\bf J}^{(1)}$,...,${\bf J}^{(k)}$ of ${\bf J}$, %such that
%and
disjoint (not necessarily nonempty) subsets ${\bf
N}^{(1)}$,...,${\bf N}^{(k)}$ of ${\bf J}$ with
\begin{equation}\label{e.Jniepus}
\Big(\bigcup_{l=1}^k {\bf J}^{(l)} \Big) \cap \Big(\bigcup_{l=1}^k
{\bf N}^{(l)} \Big) = \emptyset, \qquad \qquad  \bigcup_{l=1}^k
({\bf J}^{(l)} \cup {\bf N}^{(l)}) \neq {\bf J},
\end{equation}
and nonempty, disjoint subsets ${\bf I}^{(1)}$, ...,${\bf I}^{(k)}$ of ${\bf I}$ with $\bigcup_{l=1}^k {\bf I}^{(l)} \neq {\bf I}$ such that for $l=1,..,k$,
%$\sum_{i\in G_ j} h_i(F_i(t)) =  t$ for $j\in {\bf J}^{(l)}$ and
\begin{eqnarray}
{\bf N}^{(l)} = \{ j\in {\bf J}: G_j \subseteq \bigcup_{p=1}^l {\bf
I}^{(p)}\} \setminus (\bigcup_{p=1}^l {\bf J}^{(p)} \cup
\bigcup_{p=1}^{l-1}
{\bf N}^{(p)}), \quad \;\; \label{e.may*}\\
\sum_{i\in G_ j} h_i(F_i(t)) =  t, \qquad j\in {\bf J}^{(l)}, \qquad \qquad \nonumber \\
 \sum_{i\in G_j\cap \bigcup_{p=1}^{l}{\bf I}^{(p)}} h_i(F_i(t))   +
\sum_{i\in G_j\setminus \bigcup_{p=1}^{l} {\bf I}^{(p)}} h_i(f^{(l)})
<  t, \qquad  j\in {\bf K}^{(l)}, \label{e.neq}
\end{eqnarray}
where
\begin{equation}\label{e.frontdef}
F_i(t)=\sup\{ x\leq t:h_i(x)=h_i(f^{(p)})\}, \qquad i\in{\bf I}^{(p)}, \; p=1,...,l,
\end{equation}
and ${\bf K}^{(l)}= {\bf J}  \setminus \bigcup_{p=1}^l ({\bf
J}^{(p)} \cup {\bf N}^{(p)})$. Note that ${\bf K}^{(l)}\neq
\emptyset$ by \eqref{e.Jniepus} and $G_j\setminus \bigcup_{p=1}^{l}
{\bf I}^{(p)} \neq \emptyset$ for $j\in {\bf K}^{(l)}$, $l=1,...,k$,
by \eqref{e.may*}, so that the second sum in \eqref{e.neq} is taken
over a nonempty set of indices. (%This was done
Such numbers and sets were defined in the last paragraph
for $k=2$.) Let
%${\bf J}^{(2k)}= {\bf J}  \setminus \bigcup_{l=1}^k {\bf J}^{(2l-1)}$ and let
$f^{(k+1)}=f^{(k+1)}(t)$ be the supremum of $x\leq t$ satisfying the constraints
\begin{eqnarray}\label{e.maxk+1}
\sum_{i\in G_j\cap \bigcup_{p=1}^{k}{\bf I}^{(p)}} h_i(F_i(t))   +
\sum_{i\in G_j\setminus \bigcup_{p=1}^{k} {\bf I}^{(p)}} h_i(x)
\leq t, \qquad j \in {\bf K}^{(k)}.
\end{eqnarray}
%for every $j \in {\bf J}^{(2k)}$.
The inequality \eqref{e.neq} implies that $f^{(k+1)}>f^{(k)}$.
If $f^{(k+1)}=t$, we take
${\bf I}^{(k+1)}= {\bf I}\setminus \bigcup_{l=1}^k {\bf I}^{(l)}$, ${\bf J}^{(k+1)}={\bf K}^{(k)}$, ${\bf N}^{(k+1)}=\emptyset$ and
$F_i(t)=t$ for each
$i\in {\bf I}^{(k+1)}$,
so the definition of the vector
$F(t)$ is complete. In this case, we put $k_{max}=k+1$. If $f^{(k+1)}<t$,
%by continuity of $g^i_t(\cdot)$,
then $x=f^{(k+1)}$ satisfies \eqref{e.maxk+1}
and the set ${\bf J}^{(k+1)}$ of active constraints in \eqref{e.maxk+1} (i.e., these $j \in {\bf K}^{(k)}$, for which equality holds in \eqref{e.maxk+1} with $x=f^{(k+1)}$) is nonempty. In this case, let ${\bf I}^{(k+1)} = \bigcup_{j\in {\bf J}^{(k+1)}} G_j \setminus \bigcup_{l=1}^k {\bf I}^{(l)}$,
define ${\bf N}^{(k+1)}$ by \eqref{e.may*} with $l=k+1$
and put $F_i(t)=\sup\{ x\leq t: h_i(x)=h_i(f^{(k+1)})\}$ for $i\in{\bf I}^{(k+1)}$.
\iffalse
\begin{eqnarray*}
{\bf N}^{(k+1)} = \{ j\in {\bf J}: G_j \subseteq \bigcup_{l=1}^{k+1} {\bf
I}^{(l)}\} \setminus (\bigcup_{l=1}^{k+1} {\bf J}^{(l)} \cup
\bigcup_{l=1}^{k}
{\bf N}^{(l)}),
\\
F_i(t)=\sup\{ x\leq t: h_i(x)=h_i(f^{(k+1)})\}, \qquad i\in{\bf I}^{(k+1)}.
\end{eqnarray*}
\fi
%for $i \in {\bf I}^{(k+1)}$.
This ends the $k+1$-th step of our construction.
If $\bigcup_{l=1}^{k+1} {\bf I}^{(l)}={\bf I}$, the definition of the vector
$F(t)$ is complete. In this case, put $k_{max}=k+1$. Otherwise, we make another (i.e., the $k+2$-th) step of our algorithm, taking $k+1$ instead of $k$ and proceeding as above.

When the construction terminates after $k_{max}$ steps, we have defined the vector $F(t)$.

\begin{remark}\label{r.timedep}
{\rm The index $k_{max}$
%, the quantities $f^{(k)}$
and the sets ${\bf I}^{(k)}$, ${\bf J}^{(k)}$, ${\bf N}^{(k)}$, ${\bf K}^{(k)}$ defined above depend on the time $t$. In what follows, when we want to stress this
dependence, we write $k_{max}(t)$, ${\bf I}^{(k)}(t)$, ${\bf J}^{(k)}(t)$, ${\bf N}^{(k)}(t)$, ${\bf K}^{(k)}(t)$,
%instead of $k_{max}$, ${\bf I}^{(k)}$, ${\bf J}^{(k)}$, ${\bf N}^{(k)}$, ${\bf K}^{(k)}$,
respectively.}
\end{remark}
\iffalse
\begin{remark}\label{r.tildeF}
{\rm
The above construction implies that for every $t\geq 0$ and $j \in {\bf J}$,
\begin{equation}\label{e.max1front}
\sum_{i\in G_j} h_i(F_i(t))\leq t.
\end{equation}
}
\end{remark}
\fi
\begin{remark}\label{r.cff1}
{\rm
Let $\bar{h}(x)= \max_{j\in {\bf J}} \sum_{i\in G_j} h_i(x)$, $x\in \R$. If $\bar{h}$ is strictly increasing in $[\min_{i\in {\bf I}} x^*_i,\infty)$, then, by definition,
\begin{equation}\label{e.cff1}
f^{(1)}(t) = \bar{h}^{-1}(t) \wedge t, \qquad \qquad t\geq 0.
\end{equation}
}
\end{remark}
\begin{remark}\label{r.uprF}
{\rm
If $i\in {\bf I}^{(k)}(t)$ for some $k\in \{1,...,k_{max}(t)\}$ and if the function $h_i$ is strictly increasing in $[x^*_i,\infty)$, then %(see
(compare \eqref{e.frontdef}),
%$F_i(t) = f^{(k)} \vee x^*_i$
\begin{equation}\label{e.bigsimp}
F_i(t) = f^{(k)}(t) \vee x^*_i.
%, \qquad\qquad  i\in {\bf I}^{(k)}(t), \;\; k=1,...,k_{max}(t).
\end{equation}

}
\end{remark}
\begin{remark}\label{r.localtraffic}
{\rm
In general, some of the sets ${\bf N}^{(k)}$, $k=1,...,k_{max}$, may be nonempty, see Example \ref{ex.Flin}, to follow. If $j\in {\bf N}^{(k)}$ for some $k$, then
$$
\sum_{i\in G_j} h_i(F_i(t)) = \sum_{l=1}^k \sum_{i\in G_j\cap {\bf I}^{(l)}} h_i(f^{(l)}(t))< t.
$$
This strict inequality %is
may be interpreted as an indication of ``unavoidable bottleneck idleness'' in the corresponding locally edge minimal fluid network - transferring higher priority fluids by other %servers
resources does not allow $j$ to use its full capacity on the time interval $[0,t]$. This phenomenon is well known in the theory of resource sharing networks %with resource sharing
and it was discussed in detail, e.g., by Gurvich and Van Mieghem \cite{GVM, GVM2}.
A mild sufficient condition for all the sets ${\bf N}^{(k)}$ to be empty is that
for each $j\in {\bf J}$, $G_j \setminus \bigcup_{j'\neq j} G_{j'} \neq \emptyset$.
This corresponds to the so-called local traffic condition for the underlying network topology, under which every resource has at least one route using only that resource, see %Kang et al.
\cite{hmsy, kklw}. %, Assumption 5.1,
%or Harrison et al. \cite{hmsy}. %, p. 526.
The latter requirement is satisfied, for example, by linear networks, for which $I=J+1$ and $G_j=\{1,j+1\}$, $j=1,...,J$.
There are, however, some important %czy nie oslabic tego sformulowania?
%networks
%topologies
systems that do not satisfy the local traffic assumption, for example ring networks, %with cyclic graphs,
%sometimes
used as counterexamples for stability of the %Longest Queue First discipline
LQF protocol %(see %Birand et al.
\cite{Bis, %Dimakis and Walrand \cite{
DW}.
}
\end{remark}
% Jeszcze uwaga o alternatywnym opisie przez f^{(1)},...,f^{(J)}???

\subsection{Examples}

In this subsection, we provide two examples illustrating the construction of the mapping $F:[0,\infty)\rightarrow \R^I$ defined in Subsection \ref{ss.constr}. The first one has %very
relatively
simple structure, yielding time-independent $k_{max}$, ${\bf J}^{(k)}$, ${\bf I}^{(k)}$, ${\bf N}^{(k)}$ in $(0,\infty)$ and linear function $F$. The second one, in which $k_{max}$, ${\bf J}^{(k)}$, ${\bf I}^{(k)}$ vary in time and $F$ is nonlinear, indicates some of the difficulties encountered in
%the analysis of the
more general
%case.
situations.
\begin{example}\label{ex.Flin}
{\rm
Let $h_i(x)=x^+$ for $x\in \R$ and $i\in {\bf I}$. Then $x^*_i=0$ for each $i\in {\bf I}$. For $t\geq 0$ and $j\in {\bf J}$, \eqref{e.max1} takes the form
$|G_j| \; x^+\leq t$, so its maximal solution is
$f^{(1)}(t)= t/\max \{|G_j|:j\in {\bf J}\}$. Hence, for $t=0$ we have $f^{(1)}(0)=0$,
${\bf J}^{(1)}(0)={\bf J}$, ${\bf I}^{(1)}(0)={\bf I}$, ${\bf N}^{(1)}(0)=\emptyset$, $k_{\max}(0)=1$ and $F_i(0)=0$ for all $i$. In what follows, we assume that $t>0$. Then
\[ %begin{eqnarray*}
{\bf J}^{(1)}(t) = \{ j\in {\bf J}: |G_j| = \max \{|G_{j'}|:j'\in {\bf J} \} \},
\] %end{eqnarray*}
and $F_i(t) = \frac{t}{\max \{|G_j|:j\in {\bf J}\}}$ %, \qquad
for $i\in {\bf I}^{(1)}(t)=\bigcup_{j\in {\bf J}^{(1)}(t)} G_j$.

In the remainder of this example, the time argument in ${\bf J}^{(k)}(t)$, ${\bf I}^{(k)}(t)$, ${\bf N}^{(k)}(t)$, ${\bf K}^{(k)}(t)$, will be skipped.
If $|G_j|=|G_{j'}|$ for every $j,j'\in {\bf J}$,
then ${\bf J}^{(1)}={\bf J}$, ${\bf N}^{(1)}=\emptyset$ and
the definition of $F$ is complete. Otherwise, let ${\bf N}^{(1)}$ be given by \eqref{e.defn1}. For $j$ belonging to the set %\in
$$
{\bf K}^{(1)}=  \{j\in {\bf J}: |G_j| < \max_{j'\in {\bf J}} |G_{j'}|, G_j\setminus {\bf I}^{(1)} \neq \emptyset\},
$$
the constraint
\eqref{e.max2} takes the form
$
\frac{|G_j\cap  {\bf I}^{(1)}| \; t}{\max \{|G_{j'}|:j'\in {\bf J}\}} + |G_j\setminus {\bf I}^{(1)} | \; x^+\leq t,
%\quad j\in {\bf K}^{(1)}=  \{j\in {\bf J}: |G_j| < \max_{j'\in {\bf J}}  |G_{j'}|\},
$
yielding the maximal solution
\begin{equation}\label{e.f2const}
f^{(2)}(t) = \min_{j\in {\bf K}^{(1)}} \frac{1-\frac{|G_j\cap  {\bf I}^{(1)}| }{\max \{|G_{j'}|:j'\in {\bf J}\}}}{|G_j\setminus {\bf I}^{(1)} |} \; t.
\end{equation}
Moreover, ${\bf J}^{(2)}$ is the set of $j\in {\bf K}^{(1)}$ attaining the minimum %on the right-hand side of
in \eqref{e.f2const} and $F_i(t)= f^{(2)}(t)$ for $i\in {\bf I}^{(2)} = \bigcup_{j\in {\bf J}^{(2)}} G_j \setminus {\bf I}^{(1)}$. If $ {\bf I}^{(1)} \cup  {\bf I}^{(2)}=  {\bf I}$ (in particular, if the minimum %on the right-hand side of
in \eqref{e.f2const} is attained at every $j\in {\bf K}^{(1)}$), then the construction of $F$ is complete. Otherwise we proceed similarly, until we get $f^{(3)}(t)$,...,$f^{(k_{max})}(t)$, and hence all $F_i(t)$, $i\in {\bf I}$,
in the form of linear functions of $t$, with slopes depending (in an increasingly complicated way) on the sets $G_j$, describing the topology of the corresponding network.

Observe that %even
in this %relatively simple
example some sets ${\bf N}^{(k)}$ may, in general, be nonempty. The simplest such case occurs for $I=J=2$, $G_1=\{1\}$ and $G_2=\{1,2\}={\bf I}$. Then ${\bf J}^{(1)}=\{2\}$, ${\bf I}^{(1)}={\bf I}$, and hence ${\bf N}^{(1)}=\{1\}$ and $k_{max}=1$.
}
\end{example}
%\vspace {2mm}
In general, the main difficulty in analyzing the properties of the mapping $F$, %:[0,\infty)\rightarrow \R^I$,
already indicated in Remark \ref{r.timedep}, is the time-dependence of the index $k_{max}$ and the sets ${\bf I}^{(k)}$, ${\bf J}^{(k)}$, %${\bf N}^{(k)}$,
$k=1,...,k_{max}$. The following example, corresponding to a simple linear network topology, illustrates this point. Note that, by Remark \ref{r.localtraffic}, in this case the sets ${\bf N}^{(k)}$ are empty.
\begin{example}\label{ex.Fnotlin}
{\rm
Let $I=3$, $J=2$, $G_1=\{1,2\}$ and $G_2=\{1,3\}$. Furthermore, for $x\in \R$, let
$h_1(x)=(x+2)^+$, $h_2(x)=(x+1)^+$ and $h_3(x)=5x^+$, so that $x^*_1=-2$, $x^*_2=-1$, $x^*_3=0$.

Let $0\leq t \leq 1$. Then the maximal solution of \eqref{e.max1} is
$f^{(1)}(t)= t-2$
%. Moreover, we have
and we have  ${\bf J}^{(1)}(t)={\bf J}$. %, and t
Thus, ${\bf I}^{(1)}(t)={\bf I}$, $k_{max}(t)=1$, and \eqref{e.bigsimp} implies that
\begin{equation}\label{e.1case}
F_1(t) = t-2, \qquad \quad F_2(t)=-1, \qquad \quad F_3(t)=0.
\end{equation}

For $1< t < 7/2$, the maximal solution of \eqref{e.max1} is
$f^{(1)}(t)= t/2-3/2$, ${\bf J}^{(1)}(t)=\{1\}$, so ${\bf I}^{(1)}(t)=G_1=\{1,2\}$.
Furthermore, ${\bf K}^{(1)}(t)=\{2\}$ and the maximal solution of \eqref{e.max2} is $f^{(2)}(t)= (t-1)/10$, so ${\bf J}^{(2)}(t)=\{2\}$, ${\bf I}^{(2)}(t)=\{3\}$, $k_{max}(t)=2$ and \eqref{e.bigsimp} yields
\begin{equation}\label{e.72spot}
F_1(t) = F_2(t)= \frac{t-3}{2}, %- \frac{3}{2},
\qquad \quad
F_3(t)=\frac{t-1}{10}.
\end{equation}

For $t = 7/2$ we still have
$f^{(1)}(t)= t/2-3/2$, but ${\bf J}^{(1)}(t)={\bf J}$, ${\bf I}^{(1)}(t)={\bf I}$, $k_{max}(t)=1$. However, \eqref{e.72spot} still holds. Note that $F_1(7/2)=F_2(7/2)=F_3(7/2)=1/4$.

Finally, for $t>7/2$, the maximal solution of \eqref{e.max1} is
$f^{(1)}(t)= t/6-1/3$, ${\bf J}^{(1)}(t)=\{2\}$ and ${\bf I}^{(1)}(t)=G_2=\{1,3\}$. Thus, ${\bf K}^{(1)}(t)=\{1\}$ and the maximal solution of \eqref{e.max2} is $f^{(2)}(t)=5t/6-8/3$, so ${\bf J}^{(2)}(t)=\{1\}$, ${\bf I}^{(2)}(t)=\{2\}$, $k_{max}(t)=2$ and \eqref{e.bigsimp} yields
\begin{equation}\label{e.3case}
F_1(t) = F_3(t)= \frac{t}{6}- \frac{1}{3},  \qquad \quad F_2(t)= \frac{5t}{6}- \frac{8}{3}.
\end{equation}

It is easy to see that the mapping $F$ given by \eqref{e.1case}-\eqref{e.3case} is continuous and nondecreasing on $[0,\infty)$.
\iffalse
We will see that continuity of $F$ holds in general, as long as each function $h_i$ is strictly increasing in $[x^*_i,\infty)$, but global monotonicity of $F$ may fail for more complex underlying network topologies. See Sections \ref{s.cont} and \ref{ss.notglobal}, to follow.
\fi
This example also indicates that there is, in general, no hope for obtaining any global
%regularity properties of
results about the auxiliary functions $f^{(p)}$, $p>1$. Here, $f^{(2)}$ is linear and strictly increasing in $(1,7/2)$ and $(7/2,\infty)$, %(see Section \ref,
but it does not even exist in $[0,1]\cup \{7/2\}$.
}
\end{example}

\subsection{Partitions and inverses}\label{ss.constJ}

Let $\mathcal{J}$ be the set of ``ordered partitions'' of the set ${\bf J}$, i.e., finite sequences of subsets of ${\bf J}$ in the form $%\mathcal{D}=
(J_1,...,J_k)$, where $J_i\neq \emptyset$, $i=1,...,k$, $J_i\cap J_l = \emptyset$, $1\leq i <l\leq k$, and $\bigcup_{i=1}^k (J_i \cup N_i) = {\bf J}$,
$(\bigcup_{i=1}^k J_i) \cap (\bigcup_{i=1}^k N_i) =\emptyset$,
where the sets $N_1,...,N_k$
%$N_i$
are defined recursively as
\[
N_l = \{ j\in {\bf J}: G_j \subseteq \bigcup_{p=1}^l \bigcup_{j'\in J_p} G_{j'} \}\setminus (\bigcup_{p=1}^l J_p \cup
\bigcup_{p=1}^{l-1}
N_p), \qquad l=1,...k,
\]
(compare the first relation in \eqref{e.Jniepus} and \eqref{e.may*}).
Fix $T>0$ and %For $\mathcal{D}= (J_1,...,J_k)\in \mathcal{J}$,
let
\begin{equation}\label{e.defTD}
\mathcal{T}^\mathcal{D} = \{t\in [0,T]: {\bf J}^{(p)}(t) =J_p, p=1,...,k\}, \quad  \mathcal{D}= (J_1,...,J_k)\in \mathcal{J}.
\end{equation}
We label the nonempty sets of the form $\mathcal{T}^\mathcal{D}$, $\mathcal{D} \in \mathcal{J}$, as $\mathcal{T}_1,...\mathcal{T}_d$. Clearly, $\bigcup_{i=1}^d \mathcal{T}_i = [0,T]$.

For a function $g:\R \rightarrow [0,\infty)$ with $\lim_{x\rightarrow \infty} g(x)=\infty$, let $g^{-1}$ denote its (generalized) right-continuous inverse (see, e.g., %Whitt
\cite{whitt}, Section 13.6),
\[
g^{-1}(t) = \inf \{s\geq 0: g(s)>t\}, \qquad \quad  t\geq 0.
\]
If the function $g$ is nondecreasing and $g(x_0)=0$ for some $x_0\in \R$,
%small enough,
then
\[
g^{-1}(t) = \sup \{s\geq 0: g(s)\leq t\}, \qquad \quad  t\geq 0.
\]
Fix $\mathcal{D}= (J_1,...,J_k)\in \mathcal{J}$. For $t\in \mathcal{T}^\mathcal{D}$, we can write
$f^{(p)}(t)$, $p=1,...,k_{max}(t)$,  (and hence $F(t)$) in closed form,
using the inverse functions introduced above. To see this, let $t\in \mathcal{T}^\mathcal{D}$ and note that $k=k_{max}(t)$. Choose $j_1,...,j_k \in {\bf J}$ so that $j_p \in J_p = {\bf J}^{(p)}(t)$ for each $j=1,...,k$. Then
\begin{equation}\label{e.expinv}
f^{(1)}(t) = \Big( \sum_{i\in G_{j_1}} h_i \Big)^{-1}(t)  \; \wedge \; t,
\end{equation}
and for $p=2,...,k$, we have the recursive formulae
\begin{equation}\label{e.expinvp}
f^{(p)}(t) = \Big( \sum_{i\in G_{j_p} \cap I_p} h_i \Big)^{-1} \Big( t - \sum_{l=1}^{p-1} \sum_{i\in G_{j_p} \cap I_l} h_i(f^{(l)}(t)) \Big)  \; \wedge \; t,
\end{equation}
where
\begin{equation}\label{e.expinvsets}
I_1 = \bigcup_{j \in J_1} G_j = {\bf I}^{(1)}(t) ,  \quad
I_p = \bigcup_{j \in J_p} G_j  \setminus \bigcup_{l=1}^{p-1} I_l= {\bf I}^{(p)}(t), \quad p=2,...,k. \;
\end{equation}
(Actually, the minimization with $t$ in \eqref{e.expinv}-\eqref{e.expinvp} is necessary only in the %formula for
case of $f^{(k)}(t)$.)

\section{Continuity}\label{s.cont}
\setcounter{equation}{0}

In general,
%if some of the functions $h_i$
the function $F$ may have jumps, as the following one-dimensional
%simple
example indicates.
\begin{example}\label{e.jump}
{\rm
Let $I=J=1$ and let $h_1(x) = (x^+\wedge 1)+(x-2)^+$, $x\in \R$. Then $F_1(t)=f^{(1)}(t)=t/2$ for $0\leq t<2$ and $F_1(t)=f^{(1)}(t)=t/2+1$ for $t\geq 2$.
}
\end{example}
Clearly, the discontinuity of $F_1=f^{(1)}$ at $t=2$ in this example is caused by a ``flat spot'', i.e., the interval $[1,2]$ on which $h_1$ takes the constant value $2$, resulting in the jump of $h_1^{-1}$ (see \eqref{e.expinv}). In our queueing application, this corresponds to the lack of ``customer mass'' with deadlines in the interval $[1,2]$ in the system, causing the frontier to
jump over this empty interval. The following theorem assures that in the absence of such ``flat spots'', the function $F$ is actually continuous.
\begin{theorem}\label{t.ciag}
Assume that for every $i\in {\bf I}$, the function $h_i$ is strictly increasing in $[x^*_i,\infty)$. Then the mapping $F$ is continuous in $[0,\infty)$.
\end{theorem}
In the proof of this result, we will use the following elementary lemma.
%It is standard, but f
For the sake of completeness, we provide its justification.
\begin{lemma}\label{l.obv}
Let $(X,d)$ be a metric space. Let $y\in X$ and let $\{x_n\}$ be a sequence of elements of $X$ such that every subsequence $\{x_{n_k}\}$ of $\{x_n\}$ contains a further subsequence $\{x_{n_{k_l}}\}$ converging to $y$. Then $\lim_{n\rightarrow \infty} x_n=y$.
\end{lemma}
\noindent \textbf{Proof.} Suppose that the sequence $\{x_n\}$ does not converge to $y$. This means that there exist $\epsilon>0$ and a subsequence $\{x_{n_k}\}$ of $\{x_n\}$ such that $d(x_{n_k},y)\geq \epsilon$ for all $k$. However, we have assumed the existence of a subsequence $\{x_{n_{k_l}}\}$ of $\{x_{n_k}\}$ such that $\lim_{l\rightarrow \infty} d(x_{n_{k_l}},y)=0$, so we have a contradiction.
$\endproof$

\vspace{2mm}

The proof of Theorem \ref{t.ciag} is long and somewhat involved, so we only %present a few easy cases here and
sketch it here,
%the rest of the argument,
moving most of the technical details to the appendix.
It is convenient to introduce additional
%the
notation
\begin{equation}\label{e.D}
{\bf D}^{(p)} = {\bf J}^{(p)} \cup {\bf N}^{(p)}, \qquad p=1,...,k_{max}.
\end{equation}
Clearly, the sets ${\bf D}^{(p)}={\bf D}^{(p)}(t)$ depend on the time $t$, see Remark \ref{r.timedep}.
\vspace{2mm}

\noindent \textbf{Sketch of the proof of Theorem \ref{t.ciag}.}
Let $t_n>0$ be such that $t_n %\uparrow
\rightarrow t_0$ as $n\rightarrow \infty$.
We may assume that $t_n\leq t_0+1$ for all $n$.
Our aim is to show that $F(t_n)\rightarrow F(t_0)$, i.e., %that
for every $i\in {\bf I}$, we have
\begin{equation}\label{e.cont}
\lim_{n\rightarrow \infty} F_i(t_n) = F_i(t_0).
\end{equation}
Without
loss of generality (passing to a subsequence if necessary) we may
assume that for every $m,n\geq 1$, we have
$k_{max}(t_m)=k_{max}(t_n)$ and ${\bf J}^{(p)}(t_m)={\bf
J}^{(p)}(t_n)$ (hence ${\bf I}^{(p)}(t_m)={\bf I}^{(p)}(t_n)$, ${\bf N}^{(p)}(t_m)={\bf N}^{(p)}(t_n)$) for
$p=1,...,k_{max}(t_m)$. Consequently, in what follows, %the remainder of the proof
we will simply write $k_{max}$, instead of $k_{max}(t_n)$, $n\geq
1$, ${\bf J}^{(p)}$ instead of ${\bf J}^{(p)}(t_n)$, $n\geq 1$,
e.t.c..

%Remark \ref{r.cff1} implies that $f^{(1)}$ is strictly increasing and continuous on $[0,\infty)$.
By definition, for $n\geq 1$ and $p=1,...,k_{max}$,
$\min_{i\in {\bf I}} x^*_i \leq f^{(p)}(t_n) \leq t_n \leq t_0+1$. Hence, by Lemma \ref{l.obv}, without loss of generality (passing to a %further
subsequence if necessary) we may assume that the sequences $\{f^{(p)}(t_n)\}$,
$p=1,...,k_{max}$, converge. Let
\begin{equation}\label{e.limitp}
f^{(p)}_\infty:=\lim_{n\rightarrow \infty} f^{(p)}(t_n), \qquad \qquad   p=1,...,k_{max}.
\end{equation}
Remark \ref{r.cff1} implies that $f^{(1)}$ is %strictly increasing and
continuous, %on $[0,\infty)$,
hence
%so %it is clear that
%Clearly,
\begin{equation}\label{e.porza}
f^{(1)}(t_0)=f^{(1)}_\infty \leq f^{(2)}_\infty \leq ... \leq f^{(k_{max})}_\infty\leq t_0.
\end{equation}

First suppose that $t_0=0$.  %for $k=1,...,k_{max}$, w
We have
%$f^{(k_{max})}(t_n)$
%Using the construction of Section \ref{ss.constr} with $t=0$, we get
%Note that
\begin{equation}\label{e.f10}
f^{(1)}(0)=\min_{i\in {\bf I}} x^*_i,
\end{equation}
${\bf J}^{(1)}(0)={\bf J}$, ${\bf I}^{(1)}(0)={\bf I}$, $k_{max}(0)=1$ and $F_i(0)=x^*_i$ for all $i\in {\bf I}$. Let $i\in {\bf I}$. Then $i\in {\bf I}^{(k)}$ for some $k\in \{1,...,k_{max}\}$. By the definition of $f^{(k)}$, $0\leq h_i(f^{(k)}(t_n))\leq t_n$, and hence, as $t_n \downarrow 0$, by \eqref{e.defxi},
%\eqref{e.bigsimp},
\begin{eqnarray*}
F_i(t_n)&=& \sup\{ x %\in [f^{(k)}(t_n),t_n]
\leq t_n:h_i(x)=h_i(f^{(k)}(t_n))\} \\
%&=& \sup\{ x\in [f^{(k)}(t_n),t_n]:h_i(x)=h_i(f^{(k)}(t_n))\} \\
&\rightarrow & \sup\{ x%\in [f^{(1)}(t_0),t_0]
\leq 0:h_i(x)=0\} \; = \; x^*_i  \; = \; F_i(t_0),
\end{eqnarray*}
so $F_i$ is continuous at $0$.
Hence, in what follows, we will assume that $t_0>0$.

We will first consider the case in which $f^{(1)}(t_n)=t_n$ for
all $n\geq 1$. Then $f^{(1)}(t_0)=t_0$ and thus, for each $i\in {\bf I}$,
$ %\begin{equation}\label{e.chickencase}
F_i(t_n)=t_n\rightarrow t_0 = F_i(t_0).
$ %\end{equation}
Similarly, if $f^{(1)}(t_0)=t_0$ then $f^{(1)}(t_n) \rightarrow
t_0$, so for each $i\in {\bf I}$, the inclusion $F_i(t_n)\in
[f^{(1)}(t_n),t_n]$ implies that $F_i(t_n)\rightarrow t_0 =
F_i(t_0)$. Therefore, we may additionally assume
\begin{equation}\label{e.f<t}
f^{(1)}(t_n)<t_n, \qquad n\geq 0.
\end{equation}

We prove \eqref{e.cont} inductively for $i\in {\bf I}^{(l)}(t_0)$, $l=1,...,k_{max}(t_0)$. In the $l$-th inductive step, we define
\begin{eqnarray}
b_l &=& \min \{p=1,...,k_{max}: {\bf D}^{(p)} \setminus
\bigcup_{k=1}^{l-1} {\bf D}^{(k)}(t_0) \neq \emptyset\}, \label{e.defbl}\\
p_l &= &\max \{p=1,...,k_{max}: {\bf D}^{(p)} \cap \ {\bf
J}^{(l)}(t_0) \neq \emptyset\}. \label{e.defpl}
\end{eqnarray}
It is easy to check that
\begin{equation}\label{e.Ilt0}
{\bf I}^{(l)}(t_0)=\bigcup_{p=b_l}^{p_l} B_p^{(l)},
\end{equation}
where $B_p^{(l)}={\bf I}^{(p)} \cap {\bf I}^{(l)}(t_0)$.
We show that there exists $\bar{p}_l\in \{b_l,...,p_l \}$ such that
\begin{equation}
%\lim_{n\rightarrow \infty} f^{(k)}(t_n)
f^{(p)}_\infty= f^{(l)}(t_0), \qquad \qquad
p=b_l,...,\bar{p}_l. %\;\;  l=1,...,m,
\label{e.fblim}
\end{equation}
%and hence,
%Thus, b
By \eqref{e.bigsimp} and \eqref{e.fblim}, for $i\in B_p^{(l)}$, $p=b_l,...,\bar{p}_l$, as $n\rightarrow \infty$, we have
\begin{equation}\label{e.tak}
F_i(t_n) = f^{(p)}(t_n)\vee x^*_i \rightarrow f^{(l)}(t_0) \vee x^*_i = F_i(t_0).
\end{equation}
We also argue that if $\bar{p}_l<p_l$, then for $i\in B_p^{(l)}$, $p=\bar{p}_l+1,...,p_l$, we have
$
f^{(l)}(t_0) < f^{(p)}_\infty \leq x^*_i,
$
so by \eqref{e.bigsimp}, as $n\rightarrow \infty$,
\begin{equation}\label{e.albotak}
F_i(t_n) = f^{(p)}(t_n)\vee x^*_i \rightarrow f^{(p)}_\infty \vee x^*_i = x^*_i = f^{(l)}(t_0) \vee x^*_i=  F_i(t_0).
\end{equation}
Finally, the equations \eqref{e.Ilt0}, \eqref{e.tak}-\eqref{e.albotak} imply
\eqref{e.cont} for every $i\in {\bf I}^{(l)}(t_0)$.
The details of the above inductive argument may be found in the Appendix.

%\section{Preliminary analysis}
\section{Monotonicity}\label{s.mono}
\setcounter{equation}{0}

In this section, we investigate %the issue of
monotonicity of the mapping $F$.
It turns out that, in general, $F$ fails to be globally
%(i.e., on the whole $[0,\infty)$)
nondecreasing, even %in the case of piecewise linear functions $h_i$,
if the functions $h_i$ %, $i\in {\bf I}$,
are piecewise linear,
see Section \ref{ss.notglobal}. However, under suitable assumptions on
%the functions
$h_i$, %$i\in {\bf I}$,
%local
monotonicity of $F$
in some neighbourhood of
$0$
%near
may be established. More precisely, our goal is to find %a number
$T>0$ such that
%Fix $0\leq t <\tilde{t}<\infty$. We want to show that
for every $0\leq t <\tilde{t}<T$ and $i\in {\bf I}$, we have
\begin{equation}\label{e.ftfttilde}
F_i(t)\leq F_i(\tilde{t}).
\end{equation}
This %task is accomplished
is done in Section \ref{ss.plin} for piecewise linear functions $h_i$, $i\in {\bf I}$.

\subsection{Local monotonicity in the linear case}\label{ss.plin}
%\setcounter{equation}{0}

%Henceforth,
In this subsection we assume that for all $i\in {\bf I}$, %and $x\in \R$,
\begin{equation}\label{e.linearhi}
h_i(x)= \rho^i (x-x^*_i)^+,  \quad \qquad x\in \R,
\end{equation}
where $\rho^i$, $i\in {\bf I}$, are given positive constants.
Without loss of generality we also assume that
\begin{equation}\label{e.xorder}
x^*_1\leq x^*_2 \leq ... \leq x^*_I.
\end{equation}
Let %$n_0=0$ and let
$m^*\in \{1,...I\}$ and %let
$0=n_{0}<1\leq
n_1<n_2<...<n_{m^*}=I$ %\in \{1,...I-1\}$
be such that
\begin{equation}\label{e.x*ord}
x^*_1=...=x^*_{n_1}<x^*_{n_1+1}=...=x^*_{n_2}<...<x^*_{n_{m^*-1}+1}=...=x^*_{n_{m^*}}.
\end{equation}
Let $y^*_k=x^*_{n_k}$, $k=1,...,m^*$,
and let $y^*_{m^*+1}=\infty$. By \eqref{e.f10} and
\eqref{e.x*ord}, we have $f^{(1)}(0)=x^*_1=y^*_1$. It follows from
Remark \ref{r.cff1} that $f^{(1)}(\cdot)$ is continuous and strictly
increasing in $[0,\infty)$.
Let $t^*_1= (f^{(1)})^{-1}(y_{2}^*)$ if $m^*\geq 2$ and $t^*_1=\infty$ otherwise.
%Let $t^*_k= (f^{(1)})^{-1}(y_{k+1}^*)$
%for $k=0,...,m^*-1$ and let $t^*_{m^*}=\infty$.
%Note that $t^*_0=0$.
%We will first consider the case of $x^*_1<0$.
Let
\begin{equation}\label{e.a1}
a^{(1)}= \max_{j\in {\bf J}} \sum_{i \in G_j \cap
\{1,...,n_1\}} \rho^i.
\end{equation}
%For $t\in [0,t^*_1)$, w
We have
\begin{eqnarray}
{\bf J}^{(1)}(t) &=& \{j\in {\bf J}: \sum_{i \in G_j \cap
\{1,...,n_1\}} \rho^i = a^{(1)}\}, \qquad t\in (0,t^*_1), \label{e.J1}\\
f^{(1)}(t) &=&  (x^*_1+ t/a^{(1)}) \wedge t = (y^*_1+ t/a^{(1)})
\wedge t, \qquad t\in [0,t^*_1), \label{e.f1t}
\end{eqnarray}
(see Remark \ref{r.cff1}). Using \eqref{e.f1t},
%and solving the equation $f^{(1)}(t^*_1)=y^*_2$ for $t^*_1$,
we get $t^*_1=a^{(1)} (y^*_2-y^*_1)$.
%if $m^*\geq 2$ and $t^*_1=\infty$ otherwise.
Note that ${\bf J}^{(1)}$ (and hence ${\bf I}^{(1)}$, ${\bf
N}^{(1)}$) is constant in $(0,t^*_1)$.

If $x^*_1=0$ and $a^{(1)}\leq 1$, then $m^*=1$ and $t^*_1=\infty$. In this case,  \eqref{e.f1t} implies that $f^{(1)}(t)=t$ for $t\geq 0$, so %Hence, 
$F_i(t)=t$ for all $t\geq 0$, $i\in {\bf I}$. % so \eqref{e.ftfttilde} holds. % trivially. Consequently, i
In the remainder of this section we assume that either $x^*_1<0$, or $a^{(1)}> 1$. %so in any case
In particular, we have
\begin{equation}\label{e.f1inpractice}
f^{(1)}(t)= x^*_1+ t/a^{(1)}=y^*_1+ t/a^{(1)}  , \qquad \qquad t \in [0,t^*_1 \wedge\bar{t}_1),
\end{equation}
where $\bar{t}_1=\infty$ if $a^{(1)}\geq 1$ and $\bar{t}_1=a^{(1)} x^*_1/(a^{(1)}-1) $ otherwise.
%Since $t^*_1\leq \bar{t}$  In the latter case, $f^{(1)}(\bar{t})=\bar{t}>0\geq y^*_2=f^{(1)}(t^*_1)$, so

In what follows, we consider only
$t\in [0, t^*_1 \wedge \bar{t}_1)$.
% Tu zostawione, jak bylo!
If ${\bf I}^{(1)}={\bf I}$, then $k_{max}=1$, and hence the numbers $f^{(k)}(t)$, $k=1,..,k_{max}$, have already been defined.
Assume that  ${\bf I}^{(1)}\neq {\bf I}$. Then
we continue our construction by induction as follows.

Assume that for some $k\geq 1$, there exist strictly positive
numbers $t^*_l$, $\bar{t}_l$, $l=1,...,k$, such that for all $0< t <
\underline{t}_k:=\min_{1\leq l \leq k} (t^*_l \wedge \bar{t}_l)$ the
sets ${\bf I}^{(l)}={\bf I}^{(l)}(t)$, ${\bf J}^{(l)}={\bf
J}^{(l)}(t)$, ${\bf N}^{(l)}={\bf N}^{(l)}(t)$, $l=1,...,k$, do not
depend on $t$ and $\bigcup_{l=1}^k {\bf I}^{(l)} \neq {\bf I}$,
i.e., $k_{max}(t)>k$.
For $l=1,...,k$, let $i^{(l)}=\min ({\bf I}  \setminus \bigcup_{p=1}^{l-1} {\bf I}^{(p)})$ %, where $\bigcup_{p=1}^{0} {\bf I}^{(p)}=\emptyset$,
and let
$m^{(l)}\in \{1,...,m^*\}$ be such that $x^*_{i^{(l)}}=y^*_{m^{(l)}}$.
By definition, $1=i^{(1)}<i^{(2)}<...< i^{(k)}\leq I$ and
$1=m^{(1)}\leq m^{(2)}\leq ...\leq m^{(k)}\leq m^*$. %Let also $n^{(0)}=0$.
Furthermore, we assume that for each $l=1,...,k$, there exists a constant $a^{(l)}>0$ such that
\begin{equation}\label{e.flinpractice}
f^{(l)}(t)= x^*_{i^{(l)}}+ t/a^{(l)}=y^*_{m^{(l)}}+ t/a^{(l)} < y^*_{m^{(l)}+1}  ,  \qquad \;\;\; t \in [0,\underline{t}_k) . %, \quad l=1,...,k.
\end{equation}
%Hence, b
(For notational convenience, for $l>1$, we have defined $f^{(l)}(0)$ in \eqref{e.flinpractice} by continuity, i.e., as $f^{(l)}(0+)=x^*_{i^{(l)}}$, although ${\bf J}^{(1)}(0)={\bf J}$ and hence $f^{(l)}(0)$ has not been defined by the algorithm %described in
from Section \ref{ss.constr}.)
Note that all the above assumptions have already been verified for $k=1$.

By %\eqref{e.frontdef},
\eqref{e.bigsimp}, \eqref{e.linearhi} and \eqref{e.flinpractice},
%and the definitions of $F_i(t)$ (see ,
for $t \in (0,\underline{t}_k)$, $j\in {\bf K}^{(k)}$,
the equation \eqref{e.neq} implies %that
\begin{equation}\label{e.?}
t> \sum_{p=1}^k \sum_{i\in G_j \cap {\bf I}^{(p)}} \rho^i (f^{(p)}(t) -x^*_i)^+
=
\sum_{p=1}^k \sum_{i\in G_j \cap {\bf I}^{(p)} \cap \{ n_{m^{(p)}-1}+1,...,n_{m^{(p)}} \}} \frac{\rho^i \; t}{a^{(p)}},
\end{equation}
while the equation \eqref{e.maxk+1} takes the form
\begin{eqnarray}
t &\geq &
\sum_{p=1}^k \sum_{i\in G_j \cap {\bf I}^{(p)}} \rho^i (f^{(p)}(t) -x^*_i)^+
+
\sum_{i\in G_j\setminus \bigcup_{p=1}^{k} {\bf I}^{(p)}} \rho^i (x-x^*_i)^+ \label{e.maxk+1mod}\\
&=&
\sum_{p=1}^k \sum_{i\in G_j \cap {\bf I}^{(p)} \cap \{ n_{m^{(p)}-1}+1,...,n_{m^{(p)}} \}} \rho^i \frac{t}{a^{(p)}}
+
\sum_{i\in G_j\setminus \bigcup_{p=1}^{k} {\bf I}^{(p)}} \rho^i (x-x^*_i)^+, \nonumber
\end{eqnarray}
which, in turn, is equivalent to
\begin{equation}
\sum_{i\in G_j\setminus \bigcup_{p=1}^{k} {\bf I}^{(p)}} \rho^i(x-x^*_i)^+ \leq
b^{(k+1)}_j \ t, %\qquad j \in {\bf K}^{(1)}, \quad
\label{e.max2k}
\end{equation}
where
\[ %begin{equation}\label{e.bkp1j}
b^{(k+1)}_j= 1- \sum_{p=1}^k \frac{1}{a^{(p)}} \sum_{i\in G_j\cap {\bf I}^{(p)}
\cap \{n_{m^{(p)}-1}+1,...,n_{m^{(p)}}\}}
 \rho^i, \qquad j \in {\bf K}^{(k)}.
\] %end{equation}
By \eqref{e.?},
$b^{(k+1)}_j>0$
%is strictly positive
for each $j \in {\bf K}^{(k)}$.
Let $i^{(k+1)}=\min ({\bf I}  \setminus \bigcup_{p=1}^{k} {\bf I}^{(p)})$, let
$m^{(k+1)}\in \{1,...,m^*\}$ be such that $x^*_{i^{(k+1)}}=y^*_{m^{(k+1)}}$
and let
\begin{eqnarray*}
a^{(k+1)}_j &=& \frac{1}{b^{(k+1)}_j}  \sum_{i\in G_j\cap
\{n_{m^{(k+1)}-1}+1,...,n_{m^{(k+1)}} \} \setminus \bigcup_{p=1}^{k} {\bf I}^{(p)}}
\rho^i, \qquad j \in {\bf K}^{(k)}, \\
a^{(k+1)} &=& \max_{j \in {\bf K}^{(k)}} a^{(k+1)}_j.
\end{eqnarray*}
Recall that %for
$f^{(k+1)}(t)$ is the supremum of $x\leq t$ satisfying the constraints \eqref{e.max2k}.
%Thus, b
By an argument similar to the one used for $f^{(1)}$ above, one may check that $f^{(k+1)}(\cdot)$ is continuous and strictly increasing in $(0,\underline{t}_k)$.
Let $t^*_{k+1}=\infty$ if either $m^{(k+1)}=m^*$,
or  $f^{(k+1)}(\underline{t}_k-)\leq y_{m^{(k+1)}+1}^*$, and
$t^*_{k+1}= (f^{(k+1)})^{-1}(y_{m^{(k+1)}+1}^*)$
otherwise.
By definition, $f^{(k+1)}(0):=f^{(k+1)}(0+)=x^*_{i^{(k+1)}}=y^*_{m^{(k+1)}}$
(see the notational remark following \eqref{e.flinpractice}),
and, more generally,
for $t\leq \underline{t}_k \wedge t^*_{k+1}$,
%as long as $f^{(2)}(t)\leq y^*_{k^{(2)}+1}$,
we have
\begin{eqnarray}
{\bf J}^{(k+1)}(t) &=& \{j\in {\bf K}^{(k)}: a^{(k+1)}_j = a^{(k+1)}\} \qquad \mbox{ if } t>0,
\nonumber %\label{e.Jk+1}
\\
f^{(k+1)}(t) &=&  (x^*_{i^{(k+1)}}+ t/a^{(k+1)}) \wedge t = (y^*_{m^{(k+1)}}+ t/a^{(k+1)})
\wedge t.  \qquad  \label{e.fk+1t}
\end{eqnarray}
Using \eqref{e.fk+1t},
%and solving the equation $f^{(1)}(t^*_1)=y^*_2$ for $t^*_1$,
we get $t^*_{k+1}=a^{(k+1)} (y^*_{m^{(k+1)}+1}-y^*_{m^{(k+1)}})$
unless $t^*_{k+1}=\infty$. Note that ${\bf J}^{(k+1)}$ (and hence
${\bf I}^{(k+1)}$, ${\bf N}^{(k+1)}$) is constant in
$(0,\underline{t}_k \wedge t^*_{k+1})$.

If $y^*_{m^{(k+1)}}=0$ and $a^{(k+1)}<1$, then $m^{(k+1)}=m^*$ and $t^*_{k+1}=\infty$.
%$m^*=x^*_{n_{k^{(2)}}}$
In this case,
\eqref{e.fk+1t} implies that $f^{(k+1)}(t)=t$ for $t \in [0, \underline{t}_k)$, so $k_{max}(t)=k+1$, $F_i(t)=t$ for such $t$ and $i\in {\bf I} \setminus \bigcup_{p=1}^{k} {\bf I}^{(p)}$. 
%Consequently, \eqref{e.ftfttilde} holds in $[0, \underline{t}_k)$ for $i\notin \bigcup_{p=1}^{k} {\bf I}^{(p)}$. Hence, i
In what follows, we assume that either $y^*_{m^{(k+1)}}<0$, or $a^{(k+1)}\geq 1$. In particular,
\begin{equation}\label{e.fk+1inpractice}
f^{(k+1)}(t)= x^*_{i^{(k+1)}}+ t/a^{(k+1)}=y^*_{m^{(k+1)}}+ t/a^{(k+1)}  ,  \qquad t \in [0,\underline{t}_k \wedge t^*_{k+1}  \wedge \bar{t}_{k+1}),
\end{equation}
where $\bar{t}_{k+1}=\infty$ if $a^{(k+1)}\geq 1$ and $\bar{t}_{k+1}=a^{(k+1)} x^*_{i^{(k+1)}}/(a^{(k+1)}-1) $ otherwise.

This ends the $k+1$-th step of our construction.
If $\bigcup_{l=1}^{k+1} {\bf I}^{(l)}={\bf I}$, then for $0< t < \underline{t}_{k+1}=\min_{1\leq l \leq k+1} (t^*_l \wedge \bar{t}_l)$, we have
$k_{max}(t)=k+1$ and
the definition of
$f^{(p)}(t)$, $p=1,...,k_{max}(t)$, is complete. Otherwise, we make another (i.e., the $k+2$-th) step of our algorithm, taking $k+1$ instead of $k$ and proceeding as above.

When the construction terminates, in the time interval $(0, \underline{t}_{k_{max}})$, we have the sets
${\bf I}^{(l)}={\bf I}^{(l)}(t)$, %${\bf J}^{(l)}={\bf J}^{(l)}(t)$,
$l=1,...,k_{max}$, constant in $t$,
\iffalse
the functions $f^{(l)}(t)$, $l=1,...k_{max}-1$, given by \eqref{e.flinpractice} %(with $k=
and either $f^{(k_{max})}(t)\equiv t$, or $f^{(k_{max})}(t)$ given by \eqref{e.flinpractice} with $l=k_{max}$. In any case,
\fi
and $f^{(l)}(t)$, $l=1,...,k_{max}$, in the form of strictly increasing linear functions. Let $i\in {\bf I}$ and let $k\in \{1,...,k_{max}\}$ be such that $i\in {\bf I}^{(k)}$. Let $0\leq t < \tilde{t} < \underline{t}_{k_{max}}$. Then, by %\eqref{e.frontdef},
\eqref{e.bigsimp},
\iffalse
\begin{eqnarray*}
F_i(t) &=&\sup\{ x\leq t:h_i(x)=h_i(f^{(p)}(t))\} \\
&\leq & \sup\{ x \leq \tilde{t}:h_i(x)=h_i(f^{(p)}(\tilde{t}))\}     %\\
%&=& \sup\{ x\in [f^{(p)}(\tilde{t}),\tilde{t}]:h_i(x)=h_i(f^{(p)}(\tilde{t}))\}
\;\; = \;\; F_i(\tilde{t}),
\end{eqnarray*}
\fi
$
F_i(t) = f^{(k)}(t) \vee x^*_i \leq f^{(k)}(\tilde{t}) \vee x^*_i = F_i(\tilde{t}),
$
and \eqref{e.ftfttilde} follows.

\begin{remark}\label{r.locallin}
{\rm
Since the above argument is local in time, it actually requires only that for each $i\in {\bf I}$,
\eqref{e.linearhi} holds in some neighborhood of $x^*_i$, i.e., for all $x < X^*_i$, where $X^*_i>x^*_i$ are given constants.
(Without loss of generality we may further assume that $X^*_{i'}=X^*_{i''}$ if $i',i''\in {\bf I}$ and $x^*_{i'}=x^*_{i''}$.)
We only have to
restrict $t$ in each step of our construction to the interval $[0,\underline{t}_k')$ with
%redefine $\underline{t}_k$ as
$\underline{t}_k'=\min_{1\leq l \leq k} (t^*_l \wedge \bar{t}_l \wedge \bar{\bar{t}}_l)$, where
$\bar{\bar{t}}_l=\infty$ if $f^{(l)}(\underline{t}_l-)\leq X^*_{i^{(l)}}$ and $\bar{\bar{t}}_l= (f^{(l)})^{-1}(X^*_{i^{(l)}})$ otherwise.
In this case, %the final conclusion
\eqref{e.ftfttilde} holds for
all $i\in {\bf I}$ and
$0\leq t<\tilde{t} < \underline{t}_{k_{max}}'$. %, instead of $t\in [0, \underline{t}_{k_{max}})$.
}
\end{remark}

\subsection{Reduction lemmas and Lipschitz continuity}\label{ss.Lip}

In this subsection we consider the case of general $h_i$, assuming only that
each $h_i$ is strictly increasing in $[x^*_i,\infty)$. %Let us f
Fix $T>0$ and
recall the sets $\mathcal{T}_1,...\mathcal{T}_d$ from Subsection \ref{ss.constJ}. The following lemma reduces the problem of establishing
monotonicity of $F$ to showing its monotonicity under the additional assumption that $k_{max}$ and the sets ${\bf J}^{(p)}$, ${\bf I}^{(p)}$, $p=1,...,k_{max}$, are constant in $t$.
\begin{lemma}\label{l.redm}
Assume that for $k=1,...,d$, the mapping $F$ is nondecreasing on $\mathcal{T}_k$. Then $F$ is nondecreasing on $[0,T]$.
\end{lemma}
\textbf{Proof.}
First note that $F$ is nondecreasing on $\overline{\mathcal{T}_k}$ for each $k=1,...,d$. Indeed, let $t$, $\tilde{t}$ be such that $t< \tilde{t}$ and $t,\tilde{t} \in \overline{\mathcal{T}_k}$ for some $k$. Then for each $n \in \N$, there exist $t_n, \tilde{t}_n \in \mathcal{T}_k$ such that $t_n< \tilde{t}_n$ and $t_n \rightarrow t$, $\tilde{t}_n \rightarrow \tilde{t}$ as $n\rightarrow \infty$. By assumption, $F(t_n)\leq F(\tilde{t}_n)$ for each $n$.
Letting $n\rightarrow \infty$,
by Theorem \ref{t.ciag} we get $F(t)\leq F(\tilde{t})$, so $F$ is indeed monotone
on $\overline{\mathcal{T}_k}$.

Let $0\leq  t< \tilde{t}\leq T$ and let $k_0$ be such that $t \in \mathcal{T}_{k_0}$. Let $t_1= \sup\{ s \leq \tilde{t}: s\in \mathcal{T}_{k_0}\}$. Then $t_1 \in \overline{\mathcal{T}_{k_0}}$, so $F(t)\leq F(t_1)$. If $t_1 = \tilde{t}$, we have \eqref{e.ftfttilde} and the proof is complete. Assume that $t_1 < \tilde{t}$. In this case, $(t_1,\tilde{t}] \cap \mathcal{T}_{k_0} = \emptyset$ by the definition of $t_1$. However, there exist $k_1 \neq k_0$ and a sequence $s_n \in \mathcal{T}_{k_1}$ such that $s_n \downarrow t_1$. Let $t_2= \sup\{ s \leq \tilde{t}: s\in \mathcal{T}_{k_1}\}$. Then
$t_1,t_2 \in \overline{\mathcal{T}_{k_1}}$, so $F(t_1)\leq F(t_2)$, and hence
$F(t)\leq F(t_2)$. If $t_2 = \tilde{t}$, the proof is complete, otherwise
$(t_2,\tilde{t}] \cap (\mathcal{T}_{k_0} \cup \mathcal{T}_{k_1}) = \emptyset$. Let $k_2 \notin \{ k_0, k_1\}$ and a sequence $s_n \in \mathcal{T}_{k_2}$ be such that $s_n \downarrow t_2$. Put $t_3= \sup\{ s \leq \tilde{t}: s\in \mathcal{T}_{k_2}\}$. As above, we have $F(t_2)\leq F(t_3)$, and hence
$F(t)\leq F(t_3)$. After a finite number $l\leq d$ of such steps we get $F(t)\leq F(t_l)$ and $t_l=\tilde{t}$, so %the proof of
\eqref{e.ftfttilde}
%is complete.
holds.
$\endproof$

\vspace{2mm}

A slight modification of the above argument yields
\begin{lemma}\label{l.redL}
Assume that for $k=1,...,d$, the mapping $F$ is Lipschitz continuous on $\mathcal{T}_k$. Then $F$ is Lipschitz on $[0,T]$.
\end{lemma}
This lemma, %in turn, may be combined
together with \eqref{e.expinv}-\eqref{e.expinvsets}, may be used to provide a simple proof of the following result.
\begin{theorem}\label{t.Lipcont}
Assume that there exist $0<c<C<\infty$ such that
%for every $i\in {\bf I}$ and every $x^*_i \leq x \leq y$,
\begin{equation}\label{e.Lipcond}
c(y-x) \leq h_i(y)-h_i(x) \leq C (y-x), \qquad x^*_i \leq x \leq y, \quad i\in {\bf I}.
\end{equation}
Then $F$ is Lipschitz continuous in $[0,\infty)$.
\end{theorem}
\textbf{Proof.}
Fix $T>0$. We will argue that $F$ is Lipschitz in $[0,T]$ (with the Lipschitz constant independent on $T$). According to Lemma \ref{l.redL}, it suffices to show
this in each $\mathcal{T}_l$, $l=1,...,d$.
Fix $\mathcal{D}= (J_1,...,J_k)\in \mathcal{J}$.
%For $t\in \mathcal{T}^\mathcal{D}$,
By \eqref{e.Lipcond}, for $0\leq x \leq y$, we have
$(y-x)/(CI) \leq g(y)-g(x) \leq (y-x)/c$,
where $g$ is any of the inverses appearing in \eqref{e.expinv}-\eqref{e.expinvp}. Hence, the function $f^{(1)}$ defined by \eqref{e.expinv} is Lipschitz in $\mathcal{T}^\mathcal{D}$. Using this fact, together with \eqref{e.expinvp}, \eqref{e.Lipcond}, and proceeding by induction, we show Lipschitz continuity of $f^{(p)}$, $p=2,...,k$, in $\mathcal{T}^\mathcal{D}$. From this, we get Lipschitz continuity of $F$ in $\mathcal{T}^\mathcal{D}$ by \eqref{e.bigsimp}.
$\endproof$

\vspace{2mm}

Example \ref{e.jump} indicates that the lower bound in \eqref{e.Lipcond} is %actually
necessary for Theorem \ref{t.Lipcont} to hold.

\subsection{Local monotonicity in the $C^1$ case}\label{ss.C1mon}

In this subsection, we assume that for each $i\in {\bf I}$, $h_i\in C^1([x^*_i,\infty))$ and
\begin{equation}\label{e.roi}
\rho^i:=h_i'(x^*_i)>0,
\end{equation}
where $h_i'(x^*_i)$ denotes the right derivative of $h_i$ at $x^*_i$.
\iffalse
Let
\begin{equation}\label{e.linearhiL}
h_i^L(x)= \rho^i (x-x^*_i)^+,  \quad \qquad x\in \R, % \quad i\in {\bf I},
\end{equation}
be the ``linear part'' of the function $h_i$, $i\in {\bf I}$, near $x^*_i$ (compare \eqref{e.linearhi}). % Czy to potrzebne?
\fi Again, with no loss of generality we may assume
\eqref{e.xorder}. Define $m^*$, $n_0,...,n_{m^*}$,
$y^*_1,...,y^*_{m^*}$ as in Section \ref{ss.plin}. Because our
concern here is local monotonicity of $F$, without loss of
generality we may assume that $h_i'(x)>0$ for all $i\in {\bf I}$ and
$x\geq x^*_i$ (compare Remark \ref{r.locallin}).

The main idea of %the following
our analysis for this case
is to %treat
consider it as a small perturbation of the linear problem considered in Section \ref{ss.plin}, with $\rho^i$ %, $i\in {\bf I}$,
given by \eqref{e.roi}. In particular, let $\underline{t}_{k_{max}}$ be as in %the previous section
in Section \ref{ss.plin} and
let $k_{max}^L$ be the constant $k_{max}(t)$, $0<t<\underline{t}_{k_{max}}$,
defined there. Furthermore for $k=1,...,k_{max}^L$, let $a^{(k)}$
% for %the approximating linear algorithm our algorithm with linearized fun
be %the constants defined in
as in
%the previous section
Section \ref{ss.plin} and let ${\bf J}^{(k)}_L$, ${\bf
I}^{(k)}_L$, ${\bf N}^{(k)}_L$, ${\bf K}^{(k)}_L$, $f^{(k)}_L$, denote the sets ${\bf J}^{(k)}(t)$, ${\bf I}^{(k)}(t)$, ${\bf N}^{(k)}(t)$, ${\bf K}^{(k)}(t)$, $t\in (0,\underline{t}_{k_{max}})$, and the function $f^{(k)}$ for the
above-mentioned linear problem, respectively.

Fix $0<T\leq \underline{t}_{k_{max}}/2$ and the partition
$\mathcal{D} = (J_1,...,J_k) \in \mathcal{J}$ such that there exists
a sequence $t_n\downarrow 0$ such that $t_n\in
\mathcal{T}^\mathcal{D}$ for all $n$. Since
$f^{(1)}(t)<...<f^{(k-1)}(t)<t$ for $t\in \mathcal{T}^\mathcal{D}$,
using \eqref{e.expinv}-\eqref{e.expinvp} and proceeding by
induction, it is easy to see that %the functions
$f^{(1)}$,..., $f^{(k-1)}$ are the truncations of $C^1$ functions
%(for simplicity, also denoted by $f^{(1)}$,..., $f^{(k-1)}$ in what follows) to
to $\mathcal{T}^\mathcal{D}$, and $f^{(k)}$ is the truncation to
$\mathcal{T}^\mathcal{D}$ of a function in the form $g(t)\wedge t$,
where $g$ is $C^1$. %For simplicity of the exposition,
Therefore, in order to prove that the functions $f^{(p)}$,
$p=1,...,k=k_{max}$, (and hence, by \eqref{e.bigsimp}, $F_i$, $i\in
{\bf I}$) are nondecreasing in an intersection of
$\mathcal{T}^\mathcal{D}$ and a neighborhood of zero, it suffices to
check that
\begin{equation}\label{e.finalgoal}
(f^{(p)})'(0):= \lim_{t\downarrow 0, t\in \mathcal{T}^\mathcal{D}}
\frac{f^{(p)}(t)-f^{(p)}(0)}{t}>0, \qquad p=1,...,k, %=k_{max}.
\end{equation}
where %we have put
$f^{(p)}(0):=f^{(p)}(0+)=\lim_{n\rightarrow \infty}f^{(p)}(t_n)$. %, as in the linear case.
Consequently, by Lemma \ref{l.redm}, if we verify
\eqref{e.finalgoal}, then the proof of monotonicity of the mapping
$F$ in a neighborhood of zero is complete. In what follows, we
consider only %time points
$t %, t_n$, etc., belonging to the set $
\in \mathcal{T}^\mathcal{D}$.

We will first consider the case of $x^*_1=y^*_1=0$ (hence $m^*=1$, $n_1=I$) and $a^{(1)}\leq 1$. %Let $j\in {\bf J}^{(1)}_L$.
Then for every $j\in {\bf J}$ and $t>0$ small enough,
\begin{equation}\label{e.strangef1}
\sum_{i\in G_j} h_i(t) \leq a^{(1)} t +o(t),
\end{equation}
with equality for $j\in {\bf J}^{(1)}_L$.
If $a^{(1)}<1$, then \eqref{e.strangef1} implies that $f^{(1)}(t)=t$, and hence $F_i(t)=t$, $i\in {\bf I}$, for $t$ small enough.
If $a^{(1)}=1$, then, by \eqref{e.strangef1}, $f^{(1)}(t)=t+o(t)$. Since $f^{(1)}(t)<f^{(p)}(t) \leq t$ for $p=2,...,k$, this implies %that
$(f^{(1)})'(0) =(f^{(2)})'(0) =...=(f^{(k)})'(0) =1$. Therefore, in the remainder of the proof we may assume that either $x^*_1<0$, or $a^{(1)}> 1$.

By \eqref{e.f10}, we have $f^{(1)}(0)= x^*_1=y^*_1$.
We claim that
\begin{eqnarray}
J_1 \subseteq {\bf J}^{(1)}_L, \qquad\qquad\qquad\qquad \qquad \qquad \label{e.c1}  \\
f^{(1)}(t)= f^{(1)}_L(t)+o(t)=x^*_1+ t/a^{(1)}+o(t)=y^*_1+ t/a^{(1)}+o(t),  \;\; \label{e.c2}
\end{eqnarray}
(compare \eqref{e.f1inpractice}).
Suppose that \eqref{e.c1} is false and let $j\in J_1\setminus {\bf J}^{(1)}_L$,
$j'\in {\bf J}^{(1)}_L$.
Then, by \eqref{e.a1}-\eqref{e.J1}, for small
%, strictly positive $t\in\mathcal{T}^\mathcal{D}$, %
$t>0$ we have
\begin{eqnarray*}
t = \sum_{i \in G_{j} \cap
\{1,...,n_1\}} h_i(f^{(1)}(t)) = \sum_{i \in G_{j} \cap
\{1,...,n_1\}} \rho^i(f^{(1)}(t)-x^*_1) + o(f^{(1)}(t)-x^*_1) \qquad \\
 <  \sum_{i \in G_{j'} \cap
\{1,...,n_1\}} \rho^i(f^{(1)}(t)-x^*_1) + o(f^{(1)}(t)-x^*_1) = \sum_{i \in G_{j'} \cap
\{1,...,n_1\}} h_i(f^{(1)}(t)) \leq t.
\end{eqnarray*}
We have obtained a contradiction, proving \eqref{e.c1}. Using \eqref{e.J1} and \eqref{e.c1}, for $j\in J_1$, for small $t>0$, we get
\[
t =  \sum_{i \in G_{j} \cap
\{1,...,n_1\}} h_i(f^{(1)}(t)) = (a^{(1)}+o(t)) (f^{(1)}(t)-x^*_1),
%+ o(f^{(1)}(t)-x^*_1),
\]
yielding \eqref{e.c2}. If $k=1$, then %the proof of
\eqref{e.c2} implies
\eqref{e.finalgoal} and our proof
is complete, otherwise we proceed by induction as follows.

Assume that for some $1\leq l<k$ there are indices $r\in \{1,..,l\wedge k^L_{max}\}$,
$0=p_0<1 \leq p_1<p_2<...<p_r= l$ such that
\begin{eqnarray}
\bigcup_{p=p_{s-1}+1}^{p_s} J_p \subseteq {\bf J}^{(s)}_L, \qquad s=1,...,r, \qquad\qquad  \label{e.c1a}  \\
\bigcup_{p=p_{s-1}+1}^{p_{s}} I_p = %\bigcup_{s=1}^{m-1}
{\bf I}^{(s)}_L, \qquad s=1,...,r-1. \qquad \;\; \label{e.c3a}
\end{eqnarray}
where the sets $I_p$ are as in \eqref{e.expinvsets}.
%Moreover, w
We also assume that for $p=p_{s-1}+1,...p_s$, $s=1,...,r$, we have
\begin{equation}
f^{(p)}(t)= f^{(s)}_L(t) + o(t) = x^*_{i^{(s)}}+ t/a^{(s)}+o(t)=y^*_{m^{(s)}}+ t/a^{(s)}+o(t). \quad \label{e.c2a}
\end{equation}
By \eqref{e.c1}-\eqref{e.c2}, the above assumptions hold for $l=1$ (hence $r=1$, $p_1=1$). Note that in this case \eqref{e.c3a}
%reduces to $\emptyset = \emptyset$.
is vacuously true.

For $p=0,...,k-1$, let $K_p=\{j'\in {\bf J}: G_j \setminus \bigcup_{q=1}^p I_q \neq \emptyset\}$. Clearly, $K_0={\bf J}$. Moreover, \eqref{e.expinvsets} implies that $K_p={\bf K}^{(p)}(t)$ for $p=1,...k-1$. % we have .
We claim that
\begin{equation}\label{e.inclu}
%{\bf K}_L^{(r)} \subseteq
K_l \subseteq K_{p_{r-1}} = {\bf K}_L^{(r-1)}.
\end{equation}
Indeed,
%by \eqref{e.c1a}, $\bigcup_{p=1}^{l } J_p \subseteq \bigcup_{s=1}^{r}{\bf J}^{(s)}_L$, so $\bigcup_{p=1}^{l } I_p \subseteq \bigcup_{s=1}^{r}{\bf I}^{(s)}_L$. Thus, $\bigcup_{p=1}^{l } {\bf D}^{(p)}(t)\subseteq \bigcup_{s=1}^{r}{\bf D}^{(s)}_L$, which is equivalent to the first inclusion in \eqref{e.inclu}.
%T
the %second
inclusion in \eqref{e.inclu} is obvious, since $p_{r-1}<l$. %Finally, b
By \eqref{e.c3a},
$\bigcup_{p=1}^{p_{r-1}} I_p = \bigcup_{s=1}^{r-1}
{\bf I}^{(s)}_L$, hence $\bigcup_{p=1}^{p_{r-1}} {\bf D}^{(p)}(t)=\bigcup_{s=1}^{r-1}{\bf D}^{(s)}_L$, yielding the equality in \eqref{e.inclu}.

By \eqref{e.expinvsets} and \eqref{e.c1a} with $s=r$, $
\bigcup_{p=p_{r-1}+1}^{l} I_p \subseteq \bigcup_{j\in {\bf
J}^{(r)}_L} G_j \subseteq \bigcup_{s=1}^{r} {\bf I}^{(s)}_L$.
%On the other hand, $I_p\capI_q = \emptyset$ for $p\neq q$, so
This, together with \eqref{e.c3a} and the fact that the sets $I_p$
%, $p=1,...l$,
are disjoint,
%$I_p\cap I_q = \emptyset$ for $p\neq q$,
yields
\begin{equation}\label{e.c4a}
\bigcup_{p=p_{r-1}+1}^{l} I_p \subseteq %\bigcup_{s=1}^{m}
{\bf I}^{(r)}_L.
\end{equation}

We first assume that the inclusion in \eqref{e.c4a} is strict. Then
\begin{equation}\label{e.bomba!}
{\bf J}^{(r)}_L \cap K_l \neq \emptyset.
\end{equation}
Indeed, if \eqref{e.bomba!} is false, then ${\bf J}^{(r)}_L \subseteq \bigcup_{p=1}^{l} {\bf D}^{(p)}(t)$, and hence ${\bf I}^{(r)}_L \subseteq \bigcup_{p=1}^{l} I_p$. The latter inclusion, together with \eqref{e.c3a} and the fact that the sets $I_p$ are disjoint, yields ${\bf I}^{(r)}_L \subseteq \bigcup_{p=p_{r-1}+1}^{l} I_p$, contrary to the case assumption, so \eqref{e.bomba!} follows.
By
\eqref{e.defxi}, \eqref{e.frontdef}-\eqref{e.maxk+1}, \eqref{e.roi},
\eqref{e.c3a}-\eqref{e.c2a} and \eqref{e.c4a}, for small $t>0$, %the number
$f^{(l+1)}(t)$ is the supremum of $x\leq t$ satisfying the
constraints
\begin{eqnarray}
t &\geq& \sum_{p=1}^l \sum_{i\in G_j \cap I_p} [\rho^i
(f^{(p)}(t)-x^*_i) + o((f^{(p)}(t)-x^*_i)^+)] \nonumber\\
&&+ \sum_{i\in G_j \setminus \bigcup_{p=1}^l I_p} [\rho^i
(x-x^*_i)^+ + o((x-x^*_i)^+)] \nonumber\\
&=& \sum_{s=1}^{r-1} \sum_{i\in G_j \cap {\bf I}^{(s)}_L} [\rho^i
(f^{(s)}_L(t) -x^*_i+ o(t) ) + o((f^{(s)}_L(t) -x^*_i+ o(t) )^+)] \nonumber\\
&& + \sum_{i\in G_j \cap \bigcup_{p=p_{r-1}+1}^{l} I_p} [\rho^i
(f^{(r)}_L(t) -x^*_i+ o(t) ) + o((f^{(r)}_L(t) -x^*_i+ o(t) )^+)] \nonumber\\
&&+ \sum_{i\in G_j \setminus \bigcup_{p=1}^l I_p} [\rho^i
(x-x^*_i)^+ + o((x-x^*_i)^+)] \nonumber %\\
\end{eqnarray}
\begin{eqnarray}
&=& \sum_{s=1}^{r-1} \sum_{i\in G_j \cap {\bf I}^{(s)}_L \cap \{ n_{m^{(s)}-1}+1,...,n_{m^{(s)}} \}} \frac{\rho^i t}{a^{(s)}}  \nonumber\\
&& + \sum_{i\in G_j \cap (\bigcup_{p=p_{r-1}+1}^{l} I_p) \cap \{ n_{m^{(r)}-1}+1,...,n_{m^{(r)}} \}} \frac{\rho^i t}{a^{(r)}}
 \; + \; o(t) \nonumber\\
&&+ \sum_{i\in G_j \setminus \bigcup_{p=1}^l I_p} [\rho^i
(x-x^*_i)^+ + o((x-x^*_i)^+)], \label{e.abc}
\end{eqnarray}
for  $j\in K_l$.
%{\bf K}^{(l)}(t)=\{j'\in {\bf J}: G_j \setminus \bigcup_{p=1}^l I_p \neq \emptyset\}$.
Recall that in the corresponding linear problem, for small $t>0$, %we have
\begin{eqnarray}
t &\geq& \sum_{s=1}^{r-1} \sum_{i\in G_j \cap {\bf I}^{(s)}_L} \rho^i (f^{(s)}_L(t) -x^*_i)^+
\; +
\sum_{i\in G_j\setminus \bigcup_{s=1}^{r-1} {\bf I}^{(p)}_L} \rho^i (f^{(r)}_L(t)-x^*_i)^+ \nonumber
\\
&=&
\sum_{s=1}^{r-1} \sum_{i\in G_j \cap {\bf I}^{(s)}_L \cap \{ n_{m^{(s)}-1}+1,...,n_{m^{(s)}} \}} \frac{\rho^i \; t}{a^{(s)}}
\; +
\sum_{i\in G_j\setminus \bigcup_{s=1}^{r-1} {\bf I}^{(s)}_L} \rho^i \Big(\frac{ t}{a^{(r)}}+x^*_{i^{(r)}}-x^*_i\Big)^+ \nonumber\\
%moze usunac 2 pierwsze linijki?
&=&
\sum_{s=1}^{r-1} \sum_{i\in G_j \cap {\bf I}^{(s)}_L \cap \{ n_{m^{(s)}-1}+1,...,n_{m^{(s)}} \}} \frac{\rho^i \; t}{a^{(s)}} \nonumber\\
&& + \sum_{i\in G_j \cap (\bigcup_{p=p_{r-1}+1}^{l} I_p) \cap \{ n_{m^{(r)}-1}+1,...,n_{m^{(r)}} \}} \frac{\rho^i t}{a^{(r)}}  \nonumber \\
&&
+ \sum_{i\in G_j \setminus \bigcup_{p=1}^l I_p}
\rho^i \Big(\frac{ t}{a^{(r)}}+x^*_{i^{(r)}}-x^*_i\Big)^+ \label{e.almostthere!}
\end{eqnarray}
for  $j\in {\bf K}^{(r-1)}_L$, with equality for $j\in {\bf J}^{(r)}_L$ (compare \eqref{e.flinpractice}, \eqref{e.maxk+1mod}). Note that the second equality in \eqref{e.almostthere!} follows from \eqref{e.c4a}. Comparing \eqref{e.abc} to \eqref{e.almostthere!} and using \eqref{e.inclu}, \eqref{e.bomba!}, we get
the inclusion $J_{l+1} \subseteq {\bf J}^{(r)}_L$ and
\eqref{e.c2a} for $p=l+1$, $s=r$.
This ends the inductive step in the case of strict inclusion in \eqref{e.c4a}.

It remains to analyze the case in which
\begin{equation}\label{e.c4a=}
\bigcup_{p=p_{r-1}+1}^{l} I_p = %\bigcup_{s=1}^{m}
{\bf I}^{(r)}_L.
\end{equation}
Then \eqref{e.c3a} holds for $s=1,...,r$, so
$\bigcup_{p=1}^l I_p=\bigcup_{s=1}^{r} {\bf I}^{(p)}_L$ and
$K_l={\bf K}^{(r)}_L$ (this is the equality in \eqref{e.inclu}, with $r$ in the place of $r-1$). %, and it can be justified in the same way).
Also, \eqref{e.c3a} for $s\leq r$, together with the inequality $l<k$, implies that $r<k^L_{max}$.
The counterpart of \eqref{e.abc} in this case is
\begin{eqnarray}
t &\geq& \sum_{s=1}^{r} \sum_{i\in G_j \cap {\bf I}^{(s)}_L \cap \{ n_{m^{(s)}-1}+1,...,n_{m^{(s)}} \}} \frac{\rho^i t}{a^{(s)}}  \; + \; o(t)  \nonumber\\
&&+ \sum_{i\in G_j \setminus \bigcup_{p=1}^l I_p} [\rho^i
(x-x^*_i)^+ + o((x-x^*_i)^+)], \label{e.abc2}
\end{eqnarray}
for  $j\in K_l$. %, while the counterpart
If $r+1=k^L_{max}$, $y_{m^{(r+1)}}=y_{m^*}=0$ and $a^{(r+1)}<1$, then $f^{(r+1}_L(t) =t$ for small $t>0$ and the counterpart of \eqref{e.almostthere!} is
\begin{equation}
t > \sum_{s=1}^{r} \sum_{i\in G_j \cap {\bf I}^{(s)}_L \cap \{ n_{m^{(s)}-1}+1,...,n_{m^{(s)}} \}} \frac{\rho^i \; t}{a^{(s)}}
+ \sum_{i\in G_j \setminus \bigcup_{p=1}^l I_p}
\rho^i t \label{e.almostthere!2}
\end{equation}
for  $j\in {\bf K}^{(r)}_L=K_l$. Comparing \eqref{e.abc2} to \eqref{e.almostthere!2} we see that $f^{(l+1)}(t)=t$ and $l+1=k_{max}(t)=k$, so our inductive proof of \eqref{e.finalgoal} is complete.
If %$r+1=k^L_{max}$,
$y_{m^{(r+1)}}<0$ or $a^{(r+1)}\geq 1$, then %for small $t>0$ and
the counterpart of \eqref{e.almostthere!} is
\begin{equation}
t \geq \sum_{s=1}^{r} \sum_{i\in G_j \cap {\bf I}^{(s)}_L \cap \{ n_{m^{(s)}-1}+1,...,n_{m^{(s)}} \}} \frac{\rho^i  t}{a^{(s)}}
+ \sum_{i\in G_j \setminus \bigcup_{p=1}^l I_p}
\rho^i \Big(\frac{ t}{a^{(r+1)}}+x^*_{i^{(r+1)}}-x^*_i\Big)^+, %\;\;\;
\label{e.almostthere!3}
\end{equation}
for  $j\in {\bf K}^{(r)}_L=K_l$, with equality for $j\in {\bf
J}^{(r)}_L$. Comparing \eqref{e.abc2} to \eqref{e.almostthere!3}, we
get \eqref{e.c2a} for $p=l+1$, $s=r+1$, and $J_{l+1} \subseteq {\bf
J}^{(r+1)}_L$, yielding \eqref{e.c1a} for $s=r+1$, and the proof of
the inductive step is complete.

\subsection{Lack of global monotonicity}\label{ss.notglobal}

It is not hard to prove that, without any additional assumptions on the functions $h_i$, \eqref{e.ftfttilde} holds for every $0\leq t < \tilde{t}$ and $i \in {\bf I}^{(1)}(t)$. In spite of this, the mapping $F$ is, in general, not monotone on $[0,\infty)$, even if $h_i$, $i \in {\bf I}$, are given by \eqref{e.linearhi}, as the following example shows.
\begin{example}\label{ex.Fnotmon}
{\rm
Let $I=7$, $J=4$, $G_1=\{1,3,6\}$, $G_2=\{2,4,7\}$, $G_3=\{3,4,5,6,7\}$ and $G_4=\{6,7\}$. Next, let
\begin{eqnarray}
h_1(x) &=&h_2(x) \; =\; 2(x+11)^+, \label{e.h1} \\
h_3(x)&=&h_4(x) \; =\; (x+11)^+, \nonumber\\
h_5(x)&=&x^+, \label{e.h2}\\
h_6(x)&=&h_7(x) \; =\; 2(x+10)^+, \nonumber
\end{eqnarray}
so that $x^*_1=x^*_2=x^*_3=x_4^*=-11$, $x^*_5=0$, $x^*_6=x^*_7=-10$.
One may easily check that for $0\leq t \leq 8$, %e following assertions:
${\bf J}^{(1)}(t) = \{1,2\}$, ${\bf J}^{(2)}(t)=\{3\}$,
${\bf I}^{(1)}(t) = \{1,2,3,4,6,7\}$, ${\bf I}^{(2)}(t)=\{5\}$, ${\bf N}^{(1)}(t) = \{4\}$
and ${\bf N}^{(2)}(t) =\emptyset$. Moreover,
\begin{eqnarray}
f^{(1)}(t) =\frac{t}{3}-11, \qquad \;f^{(2)}(t) =\frac{t}{3}, \qquad \;\; \;\;0\leq t\leq 3, \label{e.f1i}\\
f^{(1)}(t) =\frac{t-53}{5}, \qquad  f^{(2)}(t) =\frac{8-t}{5}, \qquad  3\leq t\leq 8, \label{e.fi2}
\end{eqnarray}
hence %and, consequently,
\begin{equation}\label{e.f5}
F_5(t)=\frac{8-t}{5}, \qquad 3\leq t\leq 8,
\end{equation}
%so
and the mapping $F$ fails to be nondecreasing.

The ``network topology'' in the above example may be somewhat
simplified, at the price of making some of the functions $h_i(x)$, $x\geq x^*_i$, nonlinear. Namely,
let $I=5$, $J=3$, $G_1=\{1,3\}$, $G_2=\{2,4\}$ and $G_3=\{3,4,5\}$. Assume \eqref{e.h1}-\eqref{e.h2} and let
$ %begin{eqnarray*}
%h_1(x) &=&h_2(x) \; =\; 2(x+11)^+, \\
h_3(x) = h_4(x)  = (x+11)^+ + 2(x+10)^+, %\\
%h_5(x)&=&x^+,
$ %end{eqnarray*}
so that $x^*_1=x^*_2=x^*_3=x_4^*=-11$, $x^*_5=0$.
%Note that the functions $h_3$, $h_4$ are nonlinear.
(Somewhat informally, this network structure has been obtained from the previous one by removing the fourth server and merging the routes $3,6$ (resp., $4,7$) into a single route $3$ (resp., $4$).
It is easy to verify that in this case for $0\leq t \leq 8$, %e following assertions:
${\bf J}^{(1)}(t) = \{1,2\}$, ${\bf J}^{(2)}(t)=\{3\}$,
${\bf I}^{(1)}(t) = \{1,2,3,4\}$, ${\bf I}^{(2)}(t)=\{5\}$ and %. Moreover,
\eqref{e.f1i}-\eqref{e.f5} still hold.
Note that this network satisfies the local traffic condition and hence
${\bf N}^{(1)}(t) = {\bf N}^{(2)}(t) =\emptyset$ for all $t$.
}

\end{example}

\section{Further research directions}\label{s.furdir}

The results obtained in this paper may be regarded as
%are
introductory in nature and there are several important issues regarding our mapping $F$ that %need
remain to be addressed.
First, one would like to relax the assumptions on regularity of $h_i$ necessary for local monotonicity of $F$. Example \ref{ex.Fnotmon} shows that
%, in general,
``kinks'' of $h_i$ may create problems in this regard, so it is not immediately clear that the monotonicity result of Section \ref{ss.C1mon} may be carried over even to the Lipschitz case. A remedy for this problem may be creating a ``differential'' version of the algorithm from Section \ref{ss.constr}, determining the derivatives of $f^{(p)}$, rather than their values, for Lipschitz $h_i$, in a way similar to our analysis for the linear case. This, however, in the absence of $C^1$ regularity of $h_i$, yields an ODE system with discontinuous right-hand side, so even establishing existence of solutions to such a system %might
may be challenging.

Another direction that appears to be important for %queueing
applications is to skip the assumption of strict monotonicity of $h_i$. As Example \ref{e.jump} indicates, this results in jumps of the corresponding process $F$. However, from the point of view of the queueing application described in Section \ref{ss.moti}, with the functions $h_i$ given by \eqref{e.ourhi}, this %does
is not necessarily %create
a problem, because some $F_i$ may %just
``jump over the flat spots'', containing no mass of the corresponding initial distributions $\overline{Z}_i(0,\cdot)$, without causing discontinuity of the resulting locally edge minimal fluid model.
Similarly, it %might
may be useful to investigate functions $h_i$ with %(at least, isolated)
upward jumps, corrresponding to distributions with atoms. This would open an avenue to %analyze the corresponding mapping $F$
using techniques similar to those developed in %Kruk \cite{lk5}
our forthcoming paper, but %directly
for %``original''
pre-limit stochastic networks, rather than
%(or together with) %the map
%the ones
for the corresponding fluid limits.

%One of them is the characterization of

Example \ref{ex.Fnotmon} shows that there is no hope for global monotonicity of $F$ in the general case. However, it is plausible that for some simple network topologies (e.g., linear or tree networks), the mapping $F$ is monotone on $[0,\infty)$. This would greatly simplify the analysis of the corresponding fluid limits, and %possibly also
aid the investigation of the pre-limit stochastic networks.

Finally, it may be interesting to replace the relation ``$\eqslantless$'' and/or the set $A_t$ in the definition of $F(t)$ by a different partial ordering and/or admissible set, and to investigate properties and possible applications of the resulting mappings.

\section{Appendix: Inductive proof of Theorem \ref{t.ciag}}
\setcounter{equation}{0}

We continue the argument starting in Section \ref{s.cont}. As we have already explained, we may assume that $t_0>0$ and that \eqref{e.limitp}, \eqref{e.f<t} hold.

\subsection{The base case: $i\in {\bf I}^{(1)}(t_0)$}

For $j\in {\bf J}^{(1)}$, we have
$\sum_{i\in G_j} h_i(f^{(1)}(t_n)) =t_n$, so %that
$\sum_{i\in G_j} h_i(f^{(1)}(t_0)) =t_0$.
%We have shown that
Consequently, ${\bf J}^{(1)} \subseteq {\bf J}^{(1)}(t_0)$, and thus ${\bf I}^{(1)} \subseteq {\bf I}^{(1)}(t_0)$.
Let $i\in {\bf I}^{(1)}$. Then $i\in {\bf I}^{(1)}(t_0)$, so, by %\eqref{e.firstF},
\eqref{e.bigsimp},
$
F_i(t_n)= f^{(1)}(t_n) \vee x^*_i \rightarrow f^{(1)}(t_0) \vee x^*_i = F_i(t_0)
$
as $n\rightarrow \infty$.
We have shown that for $i\in {\bf I}^{(1)}$, we have \eqref{e.cont}.

If ${\bf I}^{(1)} = {\bf I}^{(1)}(t_0)$ (in particular, if ${\bf
I}^{(1)} ={\bf I}$), we have \eqref{e.cont} for every $i\in {\bf
I}^{(1)}(t_0)$. Assume that ${\bf I}^{(1)} \neq  {\bf
I}^{(1)}(t_0)$. Then
\begin{equation}\label{e.added}
{\bf J}^{(1)}(t_0)\setminus {\bf D}^{(1)}\neq \emptyset,
\end{equation}
because $\bigcup_{j\in {\bf D}^{(1)}} G_j = {\bf I}^{(1)}$ (see
\eqref{e.D}), while $\bigcup_{j\in {\bf J}^{(1)}(t_0)} G_j = {\bf
I}^{(1)}(t_0)$. Let
\begin{equation}\label{e.defp1}
p_1=\max\{p=1,...,k_{max}: {\bf D}^{(p)} \cap {\bf J}^{(1)}(t_0)
\neq \emptyset\}.
\end{equation}
By \eqref{e.added}, we have $p_1\geq 2$. For $p=1,...,p_1$, let
\begin{equation}\label{e.apbp}
A_p^{(1)}={\bf D}^{(p)} \cap {\bf J}^{(1)}(t_0), \qquad \qquad
%$B_p=\bigcup_{j\in A_p} G_j$.
B_p^{(1)}={\bf I}^{(p)} \cap {\bf I}^{(1)}(t_0).
\end{equation}
Clearly,
%\iffalse
\begin{equation}\label{e.J1t0}
{\bf J}^{(1)}(t_0)=\bigcup_{p=1}^{p_1} A_p^{(1)}
\end{equation}
and %hence
$\bigcup_{p=1}^{p_1} B_p^{(1)}\subseteq {\bf I}^{(1)}(t_0)$. We will check
that %actually
%\fi
\begin{equation}\label{e.I1t0}
{\bf I}^{(1)}(t_0)=\bigcup_{p=1}^{p_1} B_p^{(1)}.
\end{equation}
Indeed, let $i\in {\bf I}^{(1)}(t_0)$. Then $i\in G_j$ for some
$j\in {\bf J}^{(1)}(t_0)$, so, by \eqref{e.J1t0}, $j\in A_p^{(1)}
\subseteq {\bf D}^{(p)} $ for some $p\in \{1,...,p_1\}$.
Consequently, for some $q\in \{1,...,p\}$, we have $i\in {\bf
I}^{(q)}$ and hence $i\in B_q^{(1)}$. Thus, ${\bf
I}^{(1)}(t_0)\subseteq \bigcup_{p=1}^{p_1} B_p^{(1)}$ and
\eqref{e.I1t0} follows. Moreover, the above argument justifies
%For future reference, we will also
%similarly
%derive
the %following
inclusion
\begin{equation}\label{e.for!}
{\bf I}^{(1)}(t_0)\subseteq \bigcup_{p=1}^{p_1} \bigcup_{j\in
A_p^{(1)}} G_j = \bigcup_{p=1}^{p_1} \bigcup_{j\in A_p^{(1)}}
\bigcup_{q=1}^p (G_j \cap {\bf I}^{(q)}),
\end{equation}
where the equality %in \eqref{e.for!}
follows %immediately
from the inclusion $A_p^{(1)} \subseteq {\bf D}^{(p)}$.
%As we have just observed,
%every $i\in {\bf I}^{(1)}(t_0)$ belongs to some $G_j$ with $j\in
%A_p$ for some $p\in \{1,...,p_1\}$, so the inclusion in
%\eqref{e.for!} holds as well.

Let $p\in
\{2,...,p_1\}$ be such that $A_{p}^{(1)}\neq \emptyset$ and
let $j\in A_{p}^{(1)}$. %On the other hand, s
Since $j\in {\bf
J}^{(1)}(t_0)$ and \eqref{e.f<t} holds,
we have
\begin{equation}\label{e.inJ1}
\sum_{i\in G_j} h_i(f^{(1)}(t_0)) = t_0.
\end{equation}
If $f^{(p)}(t_n)=t_n$ for infinitely many $n$, then,
%$p=p_1=k_{\max}$ and,
by the definition of $f^{(p)}$ (compare
\eqref{e.frontdef}-\eqref{e.maxk+1}) and the inequality $f^{(1)}<f^{(p)}$, for these $n$ we have
\begin{eqnarray}\label{e.A1alt}
\sum_{i\in G_j\cap \bigcup_{k=1}^{p-1} {\bf I}^{(k)}}  h_i(f^{(1)}(t_n)) + \sum_{i\in
G_j\cap {\bf I}^{(p)}} h_i(t_n) &\leq& \nonumber\\
\sum_{k=1}^{p-1} \sum_{i\in G_j\cap {\bf I}^{(k)}} h_i(f^{(k)}(t_n)) + \sum_{i\in G_j\cap {\bf I}^{(p)}} h_i(t_n) &\leq& t_n.   \label{e.jinJ1}
\end{eqnarray}
Letting $n\rightarrow \infty$ in \eqref{e.A1alt}, we get
\begin{equation}\label{e.A1altlim}
\sum_{i\in G_j\cap \bigcup_{k=1}^{p-1} {\bf I}^{(k)}} h_i(f^{(1)}(t_0)) + \sum_{i\in
G_j\cap {\bf I}^{(p)}} h_i(t_0) \leq t_0.
\end{equation}
The relations \eqref{e.inJ1} and \eqref{e.A1altlim}, together with
the inequalities $f^{(1)}(t_0)\leq t_0$, $t_0>0\geq \max_{i\in {\bf
I}} x^*_i$, imply that $f^{(1)}(t_0)=t_0$, which contradicts
\eqref{e.f<t}. Consequently, since $A_{p_1}^{(1)}\neq \emptyset$ by
\eqref{e.defp1}-\eqref{e.apbp},
in the remainder of our proof we may %additionally
assume that for $p\leq p_1$,
\begin{equation}\label{e.fp<t}
f^{(p)}(t_n)<t_n, %\qquad p=1,...,p_1,
\qquad n\geq 1.
\end{equation}
Since $j\in A_{p}^{(1)} \subseteq {\bf D}^{(p)} = {\bf J}^{(p)} \cup
{\bf N}^{(p)} \subseteq {\bf K}^{(p-1)} $, for every $n$ we have
\begin{equation}\label{e.A0}
\sum_{k=1}^{p-1} \sum_{i\in G_j\cap {\bf I}^{(k)}} h_i(f^{(k)}(t_n)) + \sum_{i\in
G_j\cap {\bf I}^{(p)}} h_i(f^{(p-1)}(t_n)) < t_n,
\end{equation}
(compare \eqref{e.neq}-\eqref{e.frontdef}), and %, by \eqref{e.fp<t},
\begin{equation}\label{e.A1}
\sum_{k=1}^{p-1} \sum_{i\in G_j\cap {\bf I}^{(k)}} h_i(f^{(k)}(t_n))
+ \sum_{i\in G_j\cap {\bf I}^{(p)}} h_i(f^{(p)}(t_n)) \leq t_n, % de facto, =, ale to bedzie pokazane pozniej
\end{equation}
with equality for $j\in C_{p}^{(1)} :=  {\bf J}^{(p)} \cap {\bf
J}^{(1)}(t_0) \subseteq A_{p}^{(1)}$. Note that \eqref{e.A0} and the
above-mentioned equality in \eqref{e.A1} imply that for $n\geq 1$,
$j\in C_{p}^{(1)}$,
\begin{equation}\label{e.>xi}
f^{(p)}(t_n) > f^*_{j,p}:=\min \{x^*_i: i \in G_j\cap {\bf I}^{(p)} \} . %, \qquad n\geq 1.
\end{equation}
%Let
%\begin{equation}\label{e.limitp}
%f^{(k)}_\infty:=\lim_{n\rightarrow \infty} f^{(k)}(t_n), \qquad \qquad   k=2,...,p.
%\end{equation}
By \eqref{e.limitp} and \eqref{e.>xi},
\begin{equation}\label{e.loboundf}
f^{(p)}_\infty \geq f^*_{j,p}, \qquad \qquad j\in C_{p}^{(1)}.
\end{equation}
Letting $n\rightarrow \infty$ in \eqref{e.A1} %along this subsequence
for $j\in A_{p}^{(1)}$ and using \eqref{e.limitp}, we get
\begin{equation}\label{e.A1lim}
\sum_{i\in G_j\cap {\bf I}^{(1)}} h_i(f^{(1)}(t_0)) +
\sum_{k=2}^{p-1} \sum_{i\in G_j\cap {\bf I}^{(k)}}
h_i(f^{(k)}_\infty) + \sum_{i\in G_j\cap {\bf I}^{(p)}}
h_i(f^{(p)}_\infty)  \leq t_0.
\end{equation}
This, together with \eqref{e.inJ1}, \eqref{e.porza} and monotonicity
of $h_i$, yields the equality in \eqref{e.A1lim} for each $j\in
A_{p}^{(1)}$ and the equations
\begin{equation}\label{e.fkinf}
h_i(f^{(k)}_\infty) = h_i(f^{(1)}(t_0)), \qquad i\in G_j\cap {\bf I}^{(k)}, \;\; k=2,..,p, \;\; j\in A_p^{(1)}.
\end{equation}
By \eqref{e.for!} and \eqref{e.fkinf} with $p=p_1$, we have
\begin{equation}\label{e.fkinfpower!}
h_i(f^{(k)}_\infty) = h_i(f^{(1)}(t_0)), \qquad i\in B_k^{(1)}, \;\; k=2,..,p_1.
\end{equation}
The relations \eqref{e.porza}, \eqref{e.>xi}-\eqref{e.loboundf}, \eqref{e.fkinf} (with $k=p$) and the properties of the functions $h_i$, in turn, imply that either
\begin{equation}\label{e.C1}
f^{(2)}_\infty= ... =f^{(p)}_\infty=f^{(1)}(t_0),
\end{equation}
or
\begin{equation}\label{e.C2}
f^{(1)}(t_0)<f^{(p)}_\infty= f^*_{j,p}, \qquad \qquad j\in C_{p}^{(1)},
\end{equation}
provided that $C_{p}^{(1)}\neq \emptyset$.
In the latter case, $f^*_{j,p}$ does not depend on $j\in C_{p}^{(1)}$ and hence it will be denoted by $f^*_{p}$.

%By \eqref{e.>xi}, the numbers $f^*_{j,p}$ do not depend on the subsequence of $\{t_n\}$ under consideration. In particular, t
The above argument implies that if
$ %begin{equation}\label{e.thegoodcasep}
f^{(1)}(t_0)\geq \min_{j\in C_{p}^{(1)}} f^*_{j,p} % \qquad \mbox{ for some }j\in A_{p},
$ %end{equation}
(or if $f^*_{j,p}\neq f^*_{j',p}$ for some $j,j'\in C_{p}^{(1)}$),
then
\eqref{e.C2} cannot hold, and hence \eqref{e.C1} holds.
%\begin{equation}\label{e.ourdreamp1}
%\lim_{n\rightarrow \infty} f^{(2)}(t_n)= ...= \lim_{n\rightarrow \infty} f^{(p)}(t_n) =f^{(1)}(t_0).
%\end{equation}
Similarly, if $C_{p}^{(1)}\neq \emptyset$ and
\begin{equation}\label{e.thebadcasep}
f^{(1)}(t_0)< \min_{j\in C_{p}^{(1)}} f^*_{j,p}=f^*_{p},
\end{equation}
then, by \eqref{e.loboundf} and
\eqref{e.C2},
\begin{equation}\label{e.toobad}
%\lim_{n\rightarrow \infty} f^{(p)}(t_n)
f^{(p)}_\infty= f^*_{p}.
\end{equation}
%Let us
%We stress again that the limits \eqref{e.ourdreamp1}-\eqref{e.toobad} are independent on the choice of the subsequence of $\{t_n\}$ under consideration, and hence they are taken over the entire sequence $\{t_n\}$.

%Let us summarize our %recent
%findings as follows.
Let
$ %begin{equation}\label{e.p1bar}
\bar{p}_1=\max\{p\in \{1,...,p_1\}: %\lim_{n\rightarrow \infty} f^{(p)}(t_n)
f^{(p)}_\infty =f^{(1)}(t_0)\}.
$ %end{equation}
We have shown that
\begin{equation}\label{e.ourdreamp1bar}
%\lim_{n\rightarrow \infty} f^{(k)}(t_n)
f^{(k)}_\infty=f^{(1)}(t_0), \qquad \qquad k=1,...,\bar{p}_1,
\end{equation}
and that if $\bar{p}_1<p_1$, then for every $p>\bar{p}_1$ such that $C_p^{(1)}\neq \emptyset$, the relations \eqref{e.thebadcasep}-\eqref{e.toobad} hold.
%In particular, $A_{p_1}^{(1)}\neq \emptyset$ by \eqref{e.defp1}, %the definition of $p_1$,
%so either $\bar{p}_1=p_1$, or \eqref{e.thebadcasep}-\eqref{e.toobad} hold for $p=p_1$.
%In principle, it might 
It may be the case that
$\bar{p}_1<p_1$ and
for some $\tilde{p}\in \{\bar{p}_1+1,...,p_1\}$, we have $C_{\tilde{p}}^{(1)}= \emptyset$. %Then the above argument does {\em not} assure the existence of
%$\lim_{n\rightarrow \infty} f^{(\tilde{p})}(t_n)$, taken over the entire sequence $\{t_n\}$. %, for such $\tilde{p}$.
%However, i
In this case, if $i\in B_{\tilde{p}}^{(1)}$, then
%there exists $p\in \{\tilde{p},..., p_1\}$ such that $i\in G_j\cap {\bf I}^{(\tilde{p})}$ for some $j\in A_p^{(1)}$ (see \eqref{e.for!}). Then
the relations %\eqref{e.inJ1}
\eqref{e.porza} and
%\eqref{e.fkinf}
\eqref{e.fkinfpower!}
imply that %every limit point $f^{(\tilde{p})}_\infty$ of the sequence $\{f^{(\tilde{p})}(t_n)\}$ satisfies
\begin{equation}\label{e.intheworstcase}
f^{(1)}(t_0) < f^{(\tilde{p})}_\infty \leq x_i^*, \qquad \qquad  i\in B_{\tilde{p}}^{(1)}.
\end{equation}
% (It turns out that for $p\leq \bar{p}_1$, ${\bf J}^{(p)} \subseteq {\bf
% J}^{(1)}(t_0)$, and hence $C_{p}^{(1)}={\bf J}^{(p)}\neq \emptyset$, see \eqref{e.bigstar}, to follow).

Fix  $i\in {\bf I}^{(1)}(t_0)$.
If $i\in B_{p}^{(1)}$ for some $p\leq \bar{p}_1$, then, by \eqref{e.bigsimp} and \eqref{e.ourdreamp1bar},
\[ %begin{equation}\label{e.simplei}
F_i(t_n) = f^{(p)}(t_n) \vee x^*_i \rightarrow f^{(1)}(t_0) \vee x^*_i  =  F_i(t_0), \qquad n\rightarrow \infty,
\] %end{equation}
so \eqref{e.cont} holds.
If $i\in B_{p}^{(1)}$ for some $p> \bar{p}_1$, then there exists $k\in \{p,...,p_1\}$ such that $i\in G_j\cap {\bf I}^{(p)}$ for some $j\in A_k^{(1)}$ (see \eqref{e.for!}). If $C_p^{(1)}\neq \emptyset$, then, by \eqref{e.toobad}, %the limit
$%\lim_{n\rightarrow \infty} f^{(p)}(t_n)
f^{(p)}_\infty=f^*_p$ %exists
and, by the same argument as in \eqref{e.intheworstcase}, %we have
\begin{equation}\label{e.inthebettercase}
f^{(1)}(t_0) \leq f^*_p \leq x_i^*, \qquad \qquad  i\in B_{p}^{(1)}.
\end{equation}
Thus, since $i\in {\bf I}^{(1)}(t_0)$,
by \eqref{e.defxi}, \eqref{e.bigsimp}, \eqref{e.toobad} and \eqref{e.inthebettercase},
%and the inclusion $i\in {\bf I}^{(1)}(t_0)$,
as $n\rightarrow \infty$, %we have
\begin{eqnarray}\label{e.wymconv}
F_i(t_n) = f^{(p)}(t_n) \vee x^*_i \rightarrow f^*_p \vee x^*_i = x^*_i =
f^{(1)}(t_0) \vee x^*_i = F_i(t_0),
\end{eqnarray}
and again \eqref{e.cont} holds.
Finally, if $C_p^{(1)}= \emptyset$, then reasoning as in \eqref{e.wymconv}, but using \eqref{e.intheworstcase} instead of \eqref{e.inthebettercase}, we get
$\lim_{n\rightarrow \infty} F_i (t_n)=F_i(t_0)$.
Since, by \eqref{e.I1t0}, every $i\in {\bf I}^{(1)}(t_0)$ belongs to some $B_p^{(1)}$, $p\leq p_1$, the above argument shows that \eqref{e.cont} is true for all $i\in {\bf I}^{(1)}(t_0)$.

If ${\bf D}^{(1)}(t_0)={\bf J}$, then ${\bf I}^{(1)}(t_0)={\bf I}$ and \eqref{e.cont} holds for each $i\in {\bf I}$. In what follows, we assume that
${\bf K}^{(1)}(t_0)\neq \emptyset$.
%Czy to bedzie potrzebne????
For future reference note that
\begin{equation}\label{e.bigstar}
\bigcup_{p=1}^{\bar{p}_1} {\bf J}^{(p)} \subseteq {\bf
J}^{(1)}(t_0).
\end{equation}
%where $
Indeed, let $j\in {\bf J}^{(p)}$ for some $p\leq \bar{p}_1$. By \eqref{e.fp<t}, %we have
%$f^{(1)}(t_0)<t_0$, so, by \eqref{e.ourdreamp1bar}, we have $f^{(p)}(t_n)<t_n$ for $n$ sufficiently large. For such $n$,
$ %begin{equation}\label{e.A1in goodcases}
\sum_{i\in G_j} h_i(F_i(t_n))=\sum_{k=1}^{p} \sum_{i\in G_j\cap {\bf I}^{(k)}} h_i(f^{(k)}(t_n)) = t_n,
$ %end{equation}
which, together with \eqref{e.ourdreamp1bar}, implies that
$
\sum_{i\in G_j} h_i(f^{(1)}(t_0))= \sum_{k=1}^{p} \sum_{i\in G_j\cap {\bf I}^{(k)}} h_i(f^{(1)}(t_0)) = t_0,
$
and hence $j\in {\bf J}^{(1)}(t_0)$. %The inclusion \eqref{e.bigstar} implies that
We also have
$ %begin{equation}\label{e.bigstarD}
\bigcup_{p=1}^{\bar{p}_1} {\bf D}^{(p)} \subseteq {\bf
D}^{(1)}(t_0).
$ %end{equation}
Indeed, \eqref{e.bigstar} implies %that
$
\bigcup_{p=1}^{\bar{p}_1} {\bf I}^{(p)} = \bigcup_{p=1}^{\bar{p}_1} \bigcup_{j\in  {\bf J}^{(p)}} G_j \subseteq \bigcup_{j\in  {\bf J}^{(1)}(t_0)} G_j = {\bf
I}^{(1)}(t_0).
$
Thus, for $j\in {\bf
D}^{(p)}$, $p\leq \bar{p}_1$, we have $G_j \subseteq \bigcup_{k=1}^{p} {\bf I}^{(k)} \subseteq {\bf
I}^{(1)}(t_0)$, and hence $j\in
%{\bf N}^{(1)}(t_0)\subseteq
{\bf
D}^{(1)}(t_0)$.

%\vspace{2mm}

\subsection{Inductive assumption}

%Suppose that for some
Fix $m\in \{1,...,k_{max}(t_0)-1\}$. For $l=1,...,m$, let $b_l$, $p_l$ be given by \eqref{e.defbl}-\eqref{e.defpl}.
By definition, $b_1=1$ and $b_l\leq p_l$ for each $l=1,...,m$. %Observe that
For $l=1,...,m-1$,
let $q_{l}=b_{l+1}-1$, so that by \eqref{e.defbl} we have
\begin{equation}\label{e.inql}
\bigcup_{p=1}^{q_{l}} {\bf D}^{(p)} \subseteq \bigcup_{k=1}^{l} {\bf
D}^{(k)}(t_0), \qquad \quad \bigcup_{p=1}^{q_{l}} {\bf I}^{(p)}
\subseteq \bigcup_{k=1}^{l} {\bf I}^{(k)}(t_0).
\end{equation}
For $p=b_l,...,p_l$, $l=1,...,m$, let
\begin{equation}\label{e.apbpl}
A_p^{(l)}={\bf D}^{(p)} \cap {\bf J}^{(l)}(t_0), \quad
%$B_p=\bigcup_{j\in A_p} G_j$.
B_p^{(l)}={\bf I}^{(p)} \cap {\bf I}^{(l)}(t_0), \quad
C_p^{(l)}={\bf J}^{(p)} \cap {\bf J}^{(l)}(t_0). \;\;
\end{equation}
The definitions \eqref{e.defbl}-\eqref{e.defpl} imply that
\begin{equation}\label{e.Jlt0}
{\bf J}^{(l)}(t_0)=\bigcup_{p=b_l}^{p_l} A_p^{(l)}.
\end{equation}
(see also \eqref{e.inql}).
%Clearly, %hence
%$\bigcup_{p=b_2}^{p_2} B_p^{(l)}\subseteq {\bf I}^{(l)}(t_0)$. We will
Using \eqref{e.inql} and \eqref{e.Jlt0}, %and the second inclusion in
%it is easy to
one may easily check
%that %actually
\eqref{e.Ilt0}.

Suppose that for $l=1,...,m$, the following assertions hold.

%For $l=1,...,m$, %$b_l\leq p_l$ and
There exist indices $\bar{p}_l\in \{b_l,...,p_l \}$ such that \eqref{e.fblim} holds and
\begin{equation}
%\lim_{n\rightarrow \infty} f^{(k)}(t_n)
%f^{(k)}_\infty= f^{(l)}(t_0), \qquad \qquad
%k=b_l,...,\bar{p}_l, %\;\;  l=1,...,m,
%\label{e.fblim}\\
\bigcup_{p=b_l}^{\bar{p}_l} {\bf D}^{(p)} \subseteq
\bigcup_{p=1}^{l}{\bf D}^{(p)}(t_0). \quad \qquad \qquad
\label{e.bigstarl}
\end{equation}
Note that \eqref{e.defbl} and \eqref{e.bigstarl} imply
\begin{equation}\label{e.porzabp}
1=b_1\leq \bar{p}_1<b_2 \leq \bar{p}_2 <... <b_m \leq \bar{p}_m.
\end{equation}
%Moreover, for $l=1,...,m$, there exist indices $\bar{p}_l\in \{b_l,...,p_l \}$ such that
If $\bar{p}_l<p_l$, then for every $p\in \{\bar{p}_l+1,...,p_l\}$
such that $C_p^{(l)}\neq \emptyset$, %(in particular, for $p=p_l$),
we
have
\begin{equation}
% f^{(2)}(t_0) < f^*_p, \nonumber \\
f^{(l)}(t_0) < %\lim_{n\rightarrow \infty} f^{(p)}(t_n)
f^{(p)}_\infty
=f^*_p,
\label{e.higherlimitl}
\end{equation}
where $f^*_p$ is the common value of $f^*_{j,p}$ (defined in
\eqref{e.>xi}) for $j\in C_p^{(l)}$ and, moreover,
$ %\begin{equation}\label{e.229proxyl}
f^*_p\leq x^*_i $ %, \qquad
for $i\in B_{p}^{(l)}$.
%\end{equation}
If
% we have
$C_{\tilde{p}}^{(l)}=\emptyset$ for some $\tilde{p}\in \{ \bar{p}_l+1,..., p_l\}$, then
%every limit point $f^{(\tilde{p})}_\infty$ of the sequence $\{f^{(\tilde{p})}(t_n)\}$ satisfies
\begin{equation}\label{e.228proxyl}
f^{(l)}(t_0) < f^{(\tilde{p})}_\infty\leq x^*_i, \qquad i\in B_{\tilde{p}}^{(l)}.
\end{equation}
%To unify notation, we will use the symbol $f^{(p)}_\infty$ to denote
%$\lim_{n\rightarrow \infty} f^{(p)}(t_n)$ also for
%$p=b_l,...,\bar{p}_l$ and for $p\in \{\bar{p}_l+1,...,p_l\}$
%such that $C_p^{(l)}\neq \emptyset$.
%if the latter limit exists (and hence it is the only limit point of the sequence $\{f^{(p)}(t_n)\}$.

For $j\in A_p^{(l)}$, $p=b_l,...,p_l$ and $r=1,..,l$,
\begin{equation}\label{e.keyl}
h_i(f^{(k)}_\infty)= h_i(f^{(r)}(t_0)), \qquad i\in G_j \cap
%{\bf I}^{(k)} \cap {\bf I}^{(r)}(t_0)
B^{(r)}_k, \;\;\; k=b_r,...,p_r\wedge p.
\end{equation}
For $r=1,...,m$ and $k=b_r,...,p_r$, %we have
$
B^{(r)}_k \subseteq {\bf I}^{(r)}(t_0) \subseteq \bigcup_{j \in {\bf J}^{(r)}(t_0)} G_j = \bigcup_{p=b_r}^{p_r} \bigcup_{j \in A^{(r)}_p} G_j,
$
where the equality follows from \eqref{e.Jlt0}. Hence, \eqref{e.keyl} implies
\begin{equation}\label{e.last*}
h_i(f^{(k)}_\infty)= h_i(f^{(r)}(t_0)), \qquad i\in B^{(r)}_k, \;\; k=b_r,...,p_r, \;\; r=1,...,m.
\end{equation}
Because $f^{(1)}(t_0)<...<f^{(l)}(t_0)\leq f^{(k)}_\infty$ for $k=b_l,...,p_l$,
\eqref{e.keyl} implies also that for $p=b_l,...,p_l$,
\begin{equation}\label{e.keym}
h_i(f^{(k)}_\infty)= h_i(f^{(l)}(t_0)), \qquad i\in G_j \cap {\bf
I}^{(k)}, \;\; k=b_l,...,p, \;\; j\in A_p^{(l)}.
\end{equation}
Finally, we assume that \eqref{e.cont} holds for all $i\in \bigcup_{k=1}^{m} {\bf I}^{(k)}(t_0)$.

Note that all the above assumptions have already been verified in the case of $m=1$.

%\vspace{2mm}

\subsection{Inductive step}

%Let
%\begin{equation}\label{e.defblm}
%b_{m+1} = \min \{p=1,...,k_{max}: {\bf J}^{(p)} \setminus
%\bigcup_{k=1}^{m} {\bf J}^{(k)}(t_0) \neq \emptyset\}, \qquad q_m=b_{m+1}-1.
%\end{equation}
Define $b_{m+1}$ by \eqref{e.defbl} with $l=m+1$ and let $q_m=b_{m+1}-1$.
By definition, we have \eqref{e.inql} for $l=m$,
so ${\bf K}^{(m)}(t_0) \subseteq {\bf K}^{(q_m)}$.
The inclusion
\eqref{e.bigstarl} implies that
\begin{equation}\label{e.porzapb}
b_{m+1}> \bar{p}_m.
\end{equation}
Also, for $j\in {\bf K}^{(q_m)}$ (and hence for $j\in {\bf K}^{(m)}(t_0)$), we have
\begin{eqnarray}
 \sum_{p=1}^{q_m}\sum_{i\in G_j \cap {\bf I}^{(p)}} h_i(F_i
(t_n)) + \sum_{i\in G_j \setminus \bigcup_{p=1}^{q_m}{\bf I}^{(p)}} h_i(f^{(b_{m+1})}(t_n)) &=& \nonumber  \\
\sum_{p=1}^{q_m}\sum_{i\in G_j \cap {\bf I}^{(p)}} h_i(f^{(p)}
(t_n)) + \sum_{i\in G_j \setminus \bigcup_{p=1}^{q_m}{\bf I}^{(p)}} h_i(f^{(b_{m+1})}(t_n)) & \leq& t_n.  \qquad \label{e.ajednakm+1}  %, \qquad j\in {\bf K}^{(1)}(t_0).
\end{eqnarray}
%in which
%with equality
%holds
%for $j\in {\bf J}^{(b_m)}$.

If
\begin{equation}\label{e.easycasem+1}
f^{(b_{m+1})}(t_n)=t_n,
\end{equation}
%for infinitely many $n$,
for some $n$,
then $k_{max}=b_{m+1}$,
${\bf J}^{(b_{m+1})} = {\bf K}^{(q_{m})}$, ${\bf N}^{(b_{m+1})} = \emptyset$ and %hence
$F_i(t_n)=t_n$ for all
$i\in {\bf I}^{(b_{m+1})} = {\bf I} \setminus \bigcup_{k=1}^{q_m} {\bf I}^{(k)}$. If \eqref{e.easycasem+1} holds for infinitely many $n$, then
letting $n\rightarrow \infty$ in \eqref{e.ajednakm+1} along this subsequence, by \eqref{e.inql} for $l=m$, monotonicity of $h_i$, the inequality $F_i(t_0)\leq t_0$ for all $i\in {\bf I}$,  and the inductive assumption \eqref{e.cont} for $i\in \bigcup_{k=1}^{m} {\bf I}^{(k)}(t_0)$, we get
\begin{eqnarray}
\sum_{i\in G_j \cap\bigcup_{p=1}^{m} {\bf I}^{(p)}(t_0)} h_i(F_i
(t_0)) + \sum_{i\in G_j \setminus \bigcup_{p=1}^{m}{\bf I}^{(p)}(t_0)} h_i(t_0) &\leq&    \nonumber\\
\sum_{p=1}^{q_m}\sum_{i\in G_j \cap {\bf I}^{(p)}} h_i(F_i
(t_0)) + \sum_{i\in G_j \setminus \bigcup_{p=1}^{q_m}{\bf I}^{(p)}} h_i(t_0) &\leq& t_0 \label{e.justtocompare}
\end{eqnarray}
for all $j\in {\bf K}^{(m)}(t_0)$. This, in turn, implies that $F_i(t_0)=t_0$ for all $i\in {\bf I}^{(m+1)}(t_0)= {\bf I}\setminus \bigcup_{p=1}^{m}{\bf I}^{(p)}(t_0) \subseteq {\bf I} \setminus \bigcup_{k=1}^{q_m} {\bf I}^{(k)}$, where the inclusion follows from the second inclusion in \eqref{e.inql} for $l=m$. Hence, for all $i\in {\bf I}^{(m+1)}(t_0)$ (and thus for all $i\in {\bf I}$), we have \eqref{e.cont} along a subsequence $\{t_n\}$ satisfying \eqref{e.easycasem+1}.
Consequently, in the remainder of the proof we may assume that
$ %begin{equation}\label{e.hardcasem+1}
f^{(b_{m+1})}(t_n)<t_n$ %, \qquad
for each
$
n\geq 1,
$ %end{equation}
and thus equality holds in \eqref{e.ajednakm+1}
for $j\in {\bf J}^{(b_{m+1})}$.

We claim that
\begin{eqnarray}
%\lim_{n\rightarrow \infty} f^{(b_{m+1})}(t_n)
f^{(b_{m+1})}_\infty=f^{(m+1)}(t_0), \label{e.nextdream+1}\\
%\emptyset \neq
{\bf J}^{(b_{m+1})} \setminus \bigcup_{p=1}^{m}{\bf J}^{(p)}(t_0) \subseteq {\bf J}^{(m+1)}(t_0).  \label{e.nextdream2+1}
\end{eqnarray}
%Consider a subsequence of $\{t_n\}$ (still
%denoted by $\{t_n\}$ for convenience) along which, in addition to the limits $f^{(p)}_\infty$ of the sequences $\{f^{(p)}(t_n)\}$, $p=1,...,p_m$, we have the
%existence of the limit
%$ %\begin{equation}\label{e.limitbm+1}
%f^{(b_{m+1})}_\infty := \lim_{n\rightarrow \infty} f^{(b_{m+1})}(t_n).
%$ %\end{equation}
Letting $n\rightarrow \infty$ in \eqref{e.ajednakm+1}, we get
\begin{eqnarray}
%\sum_{p=1}^{q_m}\sum_{i\in G_j \cap {\bf I}^{(p)}} h_i(f^{(p)}_\infty) + \sum_{i\in G_j \setminus \bigcup_{p=1}^{q_m}{\bf I}^{(p)}} h_i(f^{(b_{m+1})}_\infty) &=& \nonumber\\
\sum_{p=1}^{q_m}\sum_{i\in G_j \cap {\bf I}^{(p)}} h_i(F_i
(t_0)) + \sum_{i\in G_j \setminus \bigcup_{p=1}^{q_m}{\bf I}^{(p)}} h_i(f^{(b_{m+1})}_\infty) &\leq& t_0 \label{e.inthehardcase}
\end{eqnarray}
for $j\in {\bf K}^{(q_m)}$, with equality for $j\in {\bf
J}^{(b_m+1)}$ (compare \eqref{e.justtocompare}). For $k=1,...,m$, by
\eqref{e.fblim} %, \eqref{e.higherlimitl}, \eqref{e.228proxyl}
and
%the inequality $\bar{p}_m<b_{m+1}$
\eqref{e.porzapb}, we have
\begin{equation}\label{e.usefulm+1}
f^{(k)}(t_0)=f^{(b_{k})}_\infty\leq f^{(m)}(t_0)=f^{(b_{m})}_\infty
\leq f^{(b_{m+1})}_\infty,
\end{equation}
so for $i\in {\bf I}^{(k)}(t_0)$,
$h_i(F_i(t_0))=h_i(f^{(k)}(t_0)) \leq h_i(f^{(b_{m+1})}_\infty)$. This, together with \eqref{e.inql} for $l=m$ and \eqref{e.inthehardcase} implies
\begin{equation}\label{e.hardcase2}
\sum_{i\in G_j \cap\bigcup_{p=1}^{m} {\bf I}^{(p)}(t_0)} h_i(F_i
(t_0)) + \sum_{i\in G_j \setminus \bigcup_{p=1}^{m}{\bf I}^{(p)}(t_0)} h_i(f^{(b_{m+1})}_\infty) \leq t_0
\end{equation}
for $j\in {\bf K}^{(q_m)}$ (hence for $j\in {\bf K}^{(m)}(t_0)$), and thus
\begin{equation}\label{e.onesidedm+1}
f^{(b_{m+1})}_\infty \leq f^{(m+1)}(t_0).
\end{equation}
In order to prove the opposite inequality, we will first check that
\begin{equation}\label{e.**}
{\bf J}^{(b_{m+1})}\setminus \bigcup_{p=1}^{m}{\bf J}^{(p)}(t_0) \neq \emptyset.
\end {equation}
Suppose, to the converse, that
\begin{equation}\label{e.ajeto}
{\bf J}^{(b_{m+1})}\subseteq \bigcup_{p=1}^{m}{\bf J}^{(p)}(t_0) .
\end{equation}
Then
$
{\bf I}^{(b_{m+1})} \subseteq \bigcup_{j\in  {\bf J}^{(b_{m+1})}} G_j \subseteq \bigcup_{p=1}^{m} \bigcup_{j\in  {\bf J}^{(p)}(t_0) } G_j =
 \bigcup_{p=1}^{m}{\bf I}^{(p)}(t_0) ,
$
and thus, by the second inclusion in \eqref{e.inql} for $l=m$,
$\bigcup_{p=1}^{b_{m+1}}{\bf I}^{(p)} \subseteq \bigcup_{p=1}^{m}{\bf I}^{(p)}(t_0) $. Consequently, for $j\in {\bf N}^{(b_{m+1})} $,  $G_j \subseteq \bigcup_{p=1}^{m}{\bf I}^{(p)}(t_0) $,  and hence ${\bf N}^{(b_{m+1})}  \subseteq \bigcup_{p=1}^{m}{\bf D}^{(p)}(t_0)$. This, %in turn,
together with \eqref{e.ajeto}, implies that ${\bf D}^{(b_{m+1})}  \subseteq \bigcup_{p=1}^{m}{\bf D}^{(p)}(t_0)$, which contradicts the definition of $b_{m+1}$. We have proved \eqref{e.**}.

Fix $j\in {\bf J}^{(b_{m+1})}\setminus \bigcup_{p=1}^{m}{\bf J}^{(p)}(t_0)$. Then equality holds in \eqref{e.inthehardcase}. We will check that equality holds in \eqref{e.hardcase2} as well.
%Let
Suppose that $i\in G_j \cap \bigcup_{p=1}^{m} {\bf I}^{(p)}(t_0) \setminus \bigcup_{p=1}^{q_m}{\bf I}^{(p)}$. Then $i\in {\bf I}^{(k)}(t_0) \cap {\bf I}^{(r)}$ for some $k\leq m$ and $r\geq b_{m+1}$.
%Extract a further subsequence of $\{t_n\}$ (still
%denoted by $\{t_n\}$ for convenience) along which $f^{(r)}(t_n)\rightarrow f^{(r)}_\infty$.
Letting $n\rightarrow \infty$ in the equality
$h_i(F_i(t_n))=h_i(f^{(r)}(t_n))$, and using the inductive
assumption \eqref{e.cont} for $i\in \bigcup_{p=1}^{m} {\bf
I}^{(p)}(t_0)$, we get $h_i(f^{(k)}(t_0) )= h_i(F_i(t_0))=
h_i(f^{(r)}_\infty)$, while $f^{(k)}(t_0)\leq f^{(m)}(t_0) \leq
f^{(b_{m+1})}_\infty \leq f^{(r)}_\infty$ by \eqref{e.usefulm+1}.
%This is possible only if $f^{(r)}_\infty \leq x_i^*$, and h
Hence $h_i(F_i(t_0))=h_i(f^{(b_{m+1})}_\infty)$, %= 0$. C
%and %consequently
%thus
so \eqref{e.inql} for $l=m$ and the equality in
\eqref{e.inthehardcase} imply that equality holds in
\eqref{e.hardcase2} as well. Arguing as in \eqref{e.>xi}, we get $
f^{(m+1)}(t_0)>\min \{ x^*_i:i\in G_j \setminus
\bigcup_{p=1}^{m}{\bf I}^{(p)}(t_0) \}, $ so \eqref{e.onesidedm+1}
and the equality in \eqref{e.hardcase2} imply
%As in the case of $b_2=p_1+1$, the latter equality, in turn, together with \eqref{e.hohoho} and \eqref{e.onesided2}, implies
that %actually
$f^{(b_{m+1})}_\infty = f^{(m+1)}(t_0)$.
%Since the limit point $f^{(b_{m+1})}_\infty $ of the sequence $\{f^{(b_{m+1})}(t_n)\}$ is arbitrary, w
We have %actually
proved
\eqref{e.nextdream+1}. Replacing $f^{(b_{m+1})}_\infty$ by $f^{(m+1)}(t_0)$ in the equality in \eqref{e.hardcase2} we get $j\in {\bf J}^{(m+1)}(t_0)$, so  \eqref{e.nextdream2+1} %also
holds as well.

Let $i\in {\bf I}^{(b_{m+1})} \setminus \bigcup_{p=1}^{m}{\bf I}^{(p)}(t_0)$.
We will show that \eqref{e.cont} holds. % for such $i$.
Since $i\in {\bf I}^{(b_{m+1})}$, we have $i\in G_j$ for some $j\in {\bf J}^{(b_{m+1})}$. If $j\in \bigcup_{p=1}^{m}{\bf D}^{(p)}(t_0)$, then $i\in \bigcup_{p=1}^{m}{\bf I}^{(p)}(t_0)$, contrary to the choice of %assumption on
$i$, so $j\in {\bf J}^{(b_{m+1})} \setminus \bigcup_{p=1}^{m}{\bf D}^{(p)}(t_0) \subseteq {\bf J}^{(b_{m+1})}  \setminus \bigcup_{p=1}^{m}{\bf J}^{(p)}(t_0)$ and hence, by \eqref{e.nextdream2+1}, $j\in {\bf J}^{(m+1)}(t_0)$.
Then $G_j \subseteq  \bigcup_{p=1}^{m+1}{\bf I}^{(p)}(t_0)$, so $i\in {\bf I}^{(m+1)}(t_0)$. Consequently, by \eqref{e.bigsimp} and \eqref{e.nextdream+1},
\iffalse
\begin{eqnarray}
F_i(t_n)&=& \sup\{ x\leq t_n:h_i(x)=h_i(f^{(b_{m+1})}(t_n))\} \qquad \label{e.2.30m+1}\\
&\rightarrow & \sup\{ x\leq t_0:h_i(x)=h_i(f^{(m+1)}(t_0))\} \; =\; F_i(t_0), \nonumber
\end{eqnarray}
\fi
\begin{equation}\label{e.2.30m+1}
F_i(t_n) = f^{(b_{m+1})}(t_n)\vee x^*_i \rightarrow f^{(m+1)}(t_0) \vee x^*_i =
F_i(t_0),
\end{equation}
as $n\rightarrow \infty$, so \eqref{e.cont} holds.

If ${\bf I}^{(m+1)}(t_0) \subseteq {\bf I}^{(b_{m+1})}$, %then
we have \eqref{e.cont} for each $i\in {\bf I}^{(m+1)}(t_0)$. Assume that
\begin{equation}\label{e.b2<p2+1}
{\bf I}^{(m+1)}(t_0) \setminus {\bf I}^{(b_{m+1})} \neq \emptyset
\end{equation}
and
let $p_{m+1}$ be given by \eqref{e.defpl} with $l=m+1$. By \eqref{e.inql} for $l=m$ and \eqref{e.b2<p2+1}, %we have
$b_{m+1}<p_{m+1}$. For $p=b_{m+1},...,p_{m+1}$, define the sets $A^{(m+1)}_p$, $B^{(m+1)}_p$, $C^{(m+1)}_p$ by \eqref{e.apbpl} with $l=m+1$. By the definitions of $b_{m+1}$, $p_{m+1}$ and \eqref{e.inql} for $l=m$, we have \eqref{e.Jlt0} %-\eqref{e.Ilt0}
for $l=m+1$.
Clearly, %hence
$\bigcup_{p=b_{m+1}}^{p_{m+1}} B_p^{(m+1)}\subseteq {\bf I}^{(m+1)}(t_0)$. We will check
that %actually
\eqref{e.Ilt0} holds
for $l=m+1$.
Indeed, let $i\in {\bf I}^{(m+1)}(t_0)$. Then $i\in G_j$ for some
$j\in {\bf J}^{(m+1)}(t_0)$, so, by \eqref{e.Jlt0} with $l=m+1$, $j\in A_p^{(m+1)}
\subseteq {\bf D}^{(p)} $ for some $p\in \{b_{m+1},...,p_{m+1}\}$.
Consequently, we have $i\in {\bf
I}^{(q)}$ for some $q\leq p$. However, $q>q_m$ by the second inclusion in \eqref{e.inql} with $l=m$, and hence %, by \eqref{e.defpq},
$i\in B_q^{(m+1)}$ for some $q\in \{b_{m+1},...,p_{m+1}\}$. Thus, ${\bf I}^{(m+1)}(t_0)\subseteq
\bigcup_{p=b_{m+1}}^{p_{m+1}} B_{p}^{(m+1)}$ and \eqref{e.Ilt0}
with $l=m+1$ follows.

We will first consider the case in which $f^{(m+1)}(t_0)=t_0$. Then ${\bf I}^{(m+1)}(t_0) = {\bf I} \setminus \bigcup_{p=1}^{m}{\bf I}^{(p)}(t_0)$ and $F_i(t_0)=t_0$ for each $i\in {\bf I}^{(m+1)}(t_0)$. Moreover, by \eqref{e.nextdream+1}, %we have
\begin{equation}\label{e.thedreamcomestruem+1}
%\lim_{n\rightarrow \infty} f^{(b_{m+1})}(t_n)
f^{(b_{m+1})}_\infty=t_0.
\end {equation}
\iffalse
If $b_{m+1}\leq p_m$, then \eqref{e.thedreamcomestruem+1} implies that $\lim_{n\rightarrow \infty} f^{(p_m)}(t_n)=t_0$. However, since $b_{m+1}>\bar{p}_m$, then \eqref{e.higherlimitl} holds for $p=p_m$,
%(see the discussion below \eqref{e.ourdreamp1bar}),
so $0<t_0=f^*_{p_m}\leq 0$, contradiction.
%We have thus checked that
Consequently, $b_{m+1}\geq p_m+1$, and hence %\eqref{e.bigstarplus}-
%\eqref{e.bigstarplus2} holds.
$q_m\geq p_m$.
\fi
Let $i\in {\bf I}^{(m+1)}(t_0)$. Then $i\in {\bf I}^{(p)}$ for some %$p>p_1$, so
$p\geq b_{m+1}$ by the second inclusion in \eqref{e.inql} with $l=m$, and thus
$ %begin{eqnarray*}
%t_n \geq F_i(t_n)\geq
%\sup\{ x\leq t_n:h_i(x)=h_i(f^{(b_{m+1})}(t_n))\}
f^{(b_{m+1})}(t_n) \vee x^*_i \; \leq \; F_i(t_n) \; \leq \; t_n,
$
%&\rightarrow & \sup\{ x\in [f^{(2)}(t_0),t_0]:h_i(x)=h_i(f^{(2)}(t_0))\} \; =\; %F_i(t_0),
%\end{eqnarray*}
which,
by \eqref{e.thedreamcomestruem+1}, implies that $F_i(t_n) \rightarrow t_0=F_i(t_0)$,
so \eqref{e.cont} holds for all $i$. %\in {\bf I}$.

It remains to consider the case of
\begin{equation}\label{e.fm+1<t}
f^{(m+1)}(t_0)<t_0.
\end{equation}
Let $p\in
\{b_{m+1},...,p_{m+1}\}$ be such that $A_{p}^{(m+1)}\neq \emptyset$ and
let $j\in A_{p}^{(m+1)}$. %On the other hand, s
Since $j\in {\bf
J}^{(m+1)}(t_0)$ by \eqref{e.apbpl} with $l=m+1$, the relation \eqref{e.fm+1<t}
%holds, we have
implies
\begin{equation}\label{e.inJm+1}
\sum_{k=1}^{m+1} \sum_{i\in G_j \cap {\bf I}^{(k)}(t_0)} h_i(f^{(k)}(t_0))
= t_0.
\end{equation}

Suppose that $f^{(p)}(t_n)=t_n$ for infinitely many $n$, %We claim that $p>p_m$. Indeed, i
and thus $p=k_{max}$.
Then ${\bf N}^{(p)}=\emptyset$, so ${\bf J}^{(p)}={\bf D}^{(p)}$ and %hence
$A^{(l)}_p=C^{(l)}_p$ for each $l\leq m+1$.
If  %$p\leq p_l$
$p_l=k_{max}$
for some $l\leq m$, then %$p=k_{max}$, so $l=m$ and
$p_l=p\geq b_{m+1}$, and $C^{(l)}_{p_l}=A^{(l)}_{p_l}\neq \emptyset$ by \eqref{e.defpl}, \eqref{e.apbpl}.
Hence,
by \eqref{e.porzabp}, \eqref{e.porzapb}, \eqref{e.higherlimitl}, %along this subsequence,
as $n\rightarrow \infty$, we have
$0<t_0\leftarrow t_n = f^{(p_l)}(t_n) \rightarrow f^*_{p_l}\leq 0$, contradiction. Thus, $p>\max \{p_l, l\leq m\}$, so by \eqref{e.apbpl} and \eqref{e.Ilt0},
${\bf I}^{(p)} \cap \bigcup_{l=1}^{m} {\bf I}^{(l)}(t_0)= \emptyset$.
Consequently, for $j\in A^{(m+1)}_p
%=C^{(m+1)}_p
\subseteq {\bf J}^{(m+1)}(t_0)$, $G_j \cap {\bf I}^{(p)}\subseteq {\bf I}^{(m+1)}(t_0)$ and the definition of $f^{(p)}$ (compare
\eqref{e.maxk+1}), together with \eqref{e.Ilt0} for $l=1,...,m+1$, %-\eqref{e.potrz},
imply
\begin{eqnarray}
t_n &\geq& \sum_{k=1}^{p-1} \sum_{i\in G_j\cap {\bf I}^{(k)}} h_i(f^{(k)}(t_n)) + \sum_{i\in G_j\cap {\bf I}^{(p)}} h_i(t_n) \nonumber\\
%&\geq& \sum_{i\in G_j\cap {\bf I}^{(1)}(t_0)} h_i(f^{(1)}(t_n))
%+ \sum_{k=b_2}^{p-1} \sum_{i\in G_j\cap {\bf I}^{(k)}\cap {\bf I}^{(2)}(t_0)} h_i(f^{(b_2)}(t_n)) \nonumber\\
&\geq& \sum_{l=1}^m \sum_{i\in G_j\cap {\bf I}^{(l)}(t_0)} h_i(f^{(b_l)}(t_n))
+ \sum_{k=b_{m+1}}^{p-1} \sum_{i\in G_j\cap {\bf I}^{(k)}\cap {\bf I}^{(m+1)}(t_0)} h_i(f^{(b_{m+1})}(t_n)) \nonumber\\
&&+ \sum_{i\in G_j\cap {\bf I}^{(p)}} h_i(t_n) . \nonumber %\label{e.fajnem+1}
\end{eqnarray}
Letting $n\rightarrow \infty$ and using \eqref{e.fblim}, \eqref{e.nextdream+1}, we get
\[
\sum_{l=1}^m \sum_{i\in G_j\cap {\bf I}^{(l)}(t_0)} h_i(f^{(l)}(t_0))
+ \sum_{k=b_{m+1}}^{p-1} \sum_{i\in G_j\cap {\bf I}^{(k)}\cap {\bf I}^{(m+1)}(t_0)} h_i(f^{(m+1)}(t_0)) + \sum_{i\in G_j\cap {\bf I}^{(p)}} h_i(t_0) \leq t_0,
\]
which, together with \eqref{e.inJm+1}, implies that $f^{(m+1)}(t_0)=t_0$, contrary to \eqref{e.fm+1<t}.
Consequently, since $A_{p_{m+1}}^{(m+1)}\neq \emptyset$ by \eqref{e.defpl} and \eqref{e.apbpl} with $l=m+1$, in the remainder of the proof we may assume that \eqref{e.fp<t} holds for $p=b_{m+1},...,p_{m+1}$.
This, in turn, implies \eqref{e.A0}, \eqref{e.>xi} and the equality in \eqref{e.A1} for $j\in C_p^{(m+1)}$, $p\in \{b_{m+1},...,p_{m+1}\}$, by the same argument as in the case of $j\in C_p^{(1)}$, $p\leq p_1$.

The relation
\eqref{e.>xi} implies
\begin{equation}\label{e.2.17m+1}
f^{(p)}_\infty \geq f^*_{j,p}, \quad \qquad j\in C_p^{(m+1)}, \quad p=b_{m+1},...,p_{m+1}. %de facto, nawet dla j\in{\bf J}^{(p})!!!
\end{equation}
For $j\in A_p^{(m+1)}$, $p=b_{m+1},...,p_{m+1}$,
%using \eqref{e.A1} and proceeding as in \eqref{e.fajne}, we get
letting $n\rightarrow \infty$ in \eqref{e.A1},
% along the subsequence chosen above,
we get \eqref{e.A1lim}.
%(with equality for $j\in C_p^{(m+1)}$).
This, %e latter equality,
in turn, together with \eqref{e.inJm+1}, \eqref{e.porza}, \eqref{e.defbl}-\eqref{e.defpl}, %\eqref{e.defp2},
\eqref{e.fblim}, \eqref{e.nextdream+1} and monotonicity of $h_i$,
yields \eqref{e.keyl} for $j\in A_p^{(m+1)}$,
$p=b_{m+1},...,p_{m+1}$ and $r=1,..,m+1$. Consequently,
\eqref{e.keym} holds with $l=m+1$ and $p=b_{m+1},...,p_{m+1}$.

Using \eqref{e.2.17m+1} and \eqref{e.keym} (with $l$, $p$ as above) instead of \eqref{e.loboundf} and \eqref{e.fkinf}, respectively, in %complete
analogy to the development of
\eqref{e.C1}-\eqref{e.intheworstcase},
\eqref{e.inthebettercase}, we can derive the following facts.
Let
$$
\bar{p}_{m+1}=\max \{ p\in \{b_{m+1},...,p_{m+1}\}: f^{(p)}_\infty = f^{(m+1)}(t_0) \}.
$$
Then \eqref{e.fblim} holds with $l=m+1$.
%\begin{equation}\label{e.like227m+1}
%\lim_{n\rightarrow \infty} f^{(k)}(t_n)
%f^{(k)}_\infty= f^{(m+1)}(t_0),  \qquad k=b_{m+1},...,\bar{p}_{m+1}.
%\end{equation}
%and i
If $\bar{p}_{m+1}<p_{m+1}$, then for every $p\in
\{\bar{p}_{m+1}+1,...,p_{m+1}\}$ such that $C_p^{(m+1)}\neq
\emptyset$, %(in particular, for $p=p_{m+1}$),
we have \eqref{e.higherlimitl} with $l=m+1$,
%\begin{equation}
% f^{(2)}(t_0) < f^*_p, \nonumber \\
%f^{(m+1)}(t_0) < %\lim_{n\rightarrow \infty} f^{(p)}(t_n)
%f^{(p)}_\infty=f^*_p, \label{e.higherlimitm+1}
%\end{equation}
where $f^*_p$ is the common value of $f^*_{j,p}$ for $j\in
C_p^{(m+1)}$. Moreover, for such $p$,
\begin{equation}\label{e.229proxym+1}
f^*_p\leq x^*_i, \qquad i\in B_{p}^{(m+1)}.
\end{equation}
%and the limits \eqref{e.like227m+1}-\eqref{e.higherlimitm+1} (hence the index $\bar{p}_{m+1}$) are independent on the choice of the subsequence of $\{t_n\}$,
%and hence
%so they hold over the entire sequence $\{t_n\}$.
If
% we have
$C_{\tilde{p}}^{(m+1)}=\emptyset$ for some $\tilde{p}\in \{
\bar{p}_{m+1}+1,..., p_{m+1}\}$, then \eqref{e.228proxyl} holds with $l=m+1$.
%every limit point
%$f^{(\tilde{p})}_\infty$ of the sequence $\{f^{(\tilde{p})}(t_n)\}$ satisfies
%\begin{equation}\label{e.228proxym+1}
%f^{(m+1)}(t_0) < f^{(\tilde{p})}_\infty\leq x^*_i, \qquad i\in B_{\tilde{p}}^{(m+1)}.
%\end{equation}

Let $i\in {\bf I}^{(m+1)}(t_0)$. We will check that \eqref{e.cont} holds. By \eqref{e.Ilt0}
with $l=m+1$, $i\in B^{(m+1)}_p$ for some $p\in \{b_{m+1},...,p_{m+1}\}$. If $p\leq \bar{p}_{m+1}$,
then \eqref{e.cont} follows from %\eqref{e.like227m+1} 
\eqref{e.fblim} with $l=m+1$ by an argument similar to \eqref{e.2.30m+1}. If $\bar{p}_{m+1}<p\leq p_{m+1}$, then
%since $i\in G_j$ for some $j\in {\bf I}^{(2)}(t_0)$, the equation \eqref{e.J2t0} implies that $j\in A^{(2)}_q$ for some $q\in\{p,...,p_2\}$. The remainder of the proof of \eqref{e.cont} for this case
the argument is analogous to the proof of the corresponding case
$\bar{p}_1<p\leq p_1$ for $i\in {\bf I}^{(1)}(t_0)$ (see
\eqref{e.wymconv} and the paragraph surrounding it); we just use
%\eqref{e.higherlimitm+1}
\eqref{e.higherlimitl} with $l=m+1$, \eqref{e.229proxym+1},
%\eqref{e.228proxym+1} 
\eqref{e.228proxyl} with $l=m+1$, instead of \eqref{e.toobad},
\eqref{e.inthebettercase}, \eqref{e.intheworstcase}, respectively.

If $m+1=k_{max}(t_0)$, then \eqref{e.cont} holds for all $i\in {\bf
I}$ and we are done. Assume that ${\bf K}^{(m+1)}(t_0)\neq
\emptyset$. For the sake of the next inductive step, we will show
\eqref{e.bigstarl} for $l=m+1$.
%If $m+1=k_{max}(t_0)$, then this inclusion is obvious, since its right-hand side equals ${\bf J}$.
Let $j\in {\bf J}^{(p)} \setminus \bigcup_{k=1}^m {\bf
J}^{(k)}(t_0)$ for some $b_{m+1}\leq p\leq \bar{p}_{m+1}$. Since \eqref{e.fp<t} holds for %$p=b_{m+1},...,p_{m+1}$
this $p$, we have
$ %begin{eqnarray*}
t_n= \sum_{k=1}^{p} \sum_{i\in G_j\cap {\bf I}^{(k)}} h_i(f^{(k)}(t_n)).
$ %end{eqnarray*}
Letting $n\rightarrow \infty$ and using \eqref{e.apbpl}, \eqref{e.Ilt0}, \eqref{e.defbl} with $l=m+1$, %-\eqref{e.potrz},
\eqref{e.last*} and %\eqref{e.like227m+1}
\eqref{e.fblim} with $l=m+1$,  we get
\begin{eqnarray*}
t_0&=& \sum_{k=1}^{p} \sum_{i\in G_j\cap {\bf I}^{(k)}} h_i(f^{(k)}_\infty)\\
&=& \sum_{l=1}^m \sum_{k=b_l}^{p\wedge p_l} \sum_{i\in G_j\cap {\bf I}^{(k)} \cap {\bf I}^{(l)}(t_0)} h_i(f^{(k)}_\infty) +
\sum_{k=b_{m+1}}^{p} \sum_{i\in G_j\cap {\bf I}^{(k)} \setminus \bigcup_{l=1}^m {\bf
I}^{(l)}(t_0)} h_i(f^{(k)}_\infty) \\
&=& \sum_{l=1}^m \sum_{k=b_l}^{p\wedge p_l} \sum_{i\in G_j\cap B_k^{(l)}} h_i(f^{(l)}(t_0)) +
\sum_{k=b_{m+1}}^{p} \sum_{i\in G_j\cap {\bf I}^{(k)} \setminus \bigcup_{l=1}^m {\bf
I}^{(l)}(t_0)} h_i(f^{(m+1)}(t_0)) \\
&=& \sum_{l=1}^{m} \sum_{i\in G_j\cap {\bf I}^{(l)}(t_0)} h_i(f^{(l)}(t_0))
+ \sum_{i\in G_j \setminus \bigcup_{l=1}^m {\bf
I}^{(l)}(t_0)} h_i(f^{(m+1)}(t_0)),
\end{eqnarray*}
so $j\in {\bf J}^{(m+1)}(t_0)$. We have shown that
$\bigcup_{p=b_{m+1}}^{\bar{p}_{m+1}} {\bf J}^{(p)} \subseteq
\bigcup_{p=1}^{m+1}{\bf J}^{(p)}(t_0)$, so % and hence
\begin{equation}\label{e.last-1}
\bigcup_{p=b_{m+1}}^{\bar{p}_{m+1}} {\bf I}^{(p)} \subseteq \bigcup_{p=b_{m+1}}^{\bar{p}_{m+1}} \bigcup_{j\in {\bf J}^{(p)}} G_j \subseteq
\bigcup_{p=1}^{m+1} \bigcup_{j\in {\bf J}^{(p)}(t_0)} G_j = \bigcup_{p=1}^{m+1}{\bf I}^{(p)}(t_0).
\end{equation}
Let $j\in {\bf D}^{(p)}$ for some $b_{m+1}\leq p\leq \bar{p}_{m+1}$. %Then, b
By \eqref{e.inql} for $l=m$ and \eqref{e.last-1}, %we have
\[
G_j = \Big(G_j \cap \bigcup_{p=1}^{q_{m}} {\bf I}^{(p)}\Big) \cup \Big(G_j \cap \bigcup_{p=b_{m+1}}^{p} {\bf I}^{(p)}\Big) \subseteq \bigcup_{p=1}^{m}{\bf I}^{(p)}(t_0) \; \cup
 \bigcup_{p=1}^{m+1}{\bf I}^{(p)}(t_0) = \bigcup_{p=1}^{m+1}{\bf I}^{(p)}(t_0).
\]
Thus, $j\in \bigcup_{p=1}^{m+1}{\bf D}^{(p)}(t_0)$
and \eqref{e.bigstarl} with $l=m+1$ follows.

\vspace{2mm}

In summary, we have shown that either \eqref{e.cont} holds for all $i\in {\bf I}$, or all the assertions of the inductive assumption, including \eqref{e.cont} for $i\in {\bf I}^{(l)}$, are true for $l=1,...,m+1$. This ends the inductive proof of \eqref{e.cont} for each $i\in {\bf I}$.
%$\endproof$

%\vspace{3mm}

\end{document}